\documentclass[a4paper]{article}

\usepackage[english]{babel}

\usepackage[utf8x]{inputenc}
\usepackage[T1]{fontenc}

\usepackage{multirow, longtable, caption, booktabs}
\usepackage{pdflscape}

\usepackage[a4paper,top=3cm,bottom=2cm,left=2cm,right=2cm,marginparwidth=1.75cm]{geometry}

\usepackage{amsmath}
\usepackage{mathrsfs}
\usepackage{amssymb}
\usepackage{graphicx}
\usepackage{subfigure}
\usepackage[colorinlistoftodos]{todonotes}
\usepackage[colorlinks=true, allcolors=blue]{hyperref}
\usepackage{setspace}

\usepackage{ifthen}
\newboolean{Arxiv}
\setboolean{Arxiv}{true}

\usepackage{dsfont}
\usepackage{natbib}
\usepackage{nameref}
\usepackage{amsthm}
\newtheorem{theorem}{Theorem}
\newtheorem{proposition}{Proposition}
\newtheorem{lemma}{Lemma}

\newtheorem{definition}{Definition}

\newtheorem{assumption}{Assumption}

\DeclareMathOperator*{\argmin}{argmin}

\usepackage[ruled, linesnumbered]{algorithm2e}
\usepackage{threeparttable}

\usepackage{bm}

\global\long\def\inprod#1#2{\big\langle #1,#2\big\rangle }%

\global\long\def\inner#1#2{\langle#1,#2\rangle}%

\global\long\def\norm#1{\Vert#1\Vert}%
\global\long\def\bnorm#1{\big\Vert#1\big\Vert}%
\global\long\def\Bnorm#1{\Big\Vert#1\Big\Vert}%

\global\long\def\brbra#1{\big(#1\big)}%
\global\long\def\abs#1{\vert#1\vert}%
\global\long\def\bcbra#1{\big\{#1\big\}}%
\global\long\def\vertiii#1{\left\vert \kern-0.25ex  \left\vert \kern-0.25ex  \left\vert #1\right\vert \kern-0.25ex  \right\vert \kern-0.25ex  \right\vert }%
\global\long\def\mcal#1{\mathcal{#1}}%
\global\long\def\mbb#1{\mathbb{#1}}%
\global\long\def\argmin{\operatornamewithlimits{argmin}}%
\global\long\def\st{\operatornamewithlimits{s.t.}}%
\global\long\def\and{\mathrm{and}}%
\global\long\def\vec{\operatorname{vec}}%
\global\long\def\rone{\textrm{I}}%
\global\long\def\rtwo{\textrm{II}}%
\global\long\def\rthree{\textrm{III}}%

\global\long\def\BM{\text{Burer-Monteiro}}%
\global\long\def\WS{\text{warm-start}}%

\newcommand{\oursolver}{LoRADS}
\newcommand{\admm}{LRADMM}
\newcommand{\alm}{LRALM}

\title{A Low-Rank ADMM Splitting Approach for Semidefinite Programming}

\author{
\begin{tabular}{cccc}
Qiushi Han\thanks{University of Illinois Urbana-Champaign, \href{mailto:joshhan2@illinois.edu}{joshhan2@illinois.edu}}  &
Chenxi Li\thanks{Shanghai University of Finance and Economics, \{\href{mailto:chenxili@stu.sufe.edu.cn}{chenxili@stu.sufe.edu.cn}, \href{mailto:zhenweilin@163.sufe.edu.cn}{zhenweilin@163.sufe.edu.cn}\}} &
Zhenwei Lin\footnotemark[2]  &
Caihua Chen\thanks{Nanjing University, \href{mailto:chchen@nju.edu.cn}{chchen@nju.edu.cn}} \\
Qi Deng\thanks{Shanghai Jiao Tong University, \{\href{mailto:qdeng24@sjtu.edu.cn}{qdeng24@sjtu.edu.cn}, \href{mailto:ddge@sjtu.edu.cn}{ddge@sjtu.edu.cn}, \href{mailto:hkl1u@sjtu.edu.cn}{hkl1u@sjtu.edu.cn}\}} &
Dongdong Ge\footnotemark[4] &
Huikang Liu\footnotemark[4] &
Yinyu Ye\thanks{Stanford University, \href{mailto:yinyu-ye@stanford.edu}{yinyu-ye@stanford.edu}}
\end{tabular}
}

\begin{document}
\maketitle

\begin{abstract}
We introduce a new first-order method for solving general semidefinite programming problems,  based on the alternating direction method of multipliers (ADMM) and a matrix-splitting technique.
Our algorithm has an advantage over the Burer-Monteiro approach as it only involves much easier quadratically regularized subproblems in each iteration. 
For a linear objective, the subproblems are well-conditioned quadratic programs that can be efficiently solved by the standard conjugate gradient method. 
We show that the ADMM algorithm achieves sublinear or linear convergence rates to the KKT solutions under different conditions. 
Building on this theoretical development, we present LoRADS, a new solver for linear SDP based on the \underline{\bf Lo}w-\underline{\bf R}ank \underline{\bf AD}MM \underline{\bf S}plitting approach.
LoRADS incorporates several strategies that significantly increase its efficiency.
Firstly, it initiates with a warm-start phase that uses the Burer-Monteiro approach.
Moreover, motivated by the SDP low-rank theory~\citep{soandye2008},
LoRADS chooses an initial rank of logarithmic order and then employs a dynamic approach to increase the rank. 
Numerical experiments indicate that LoRADS exhibits promising performance on various SDP problems. A noteworthy achievement of LoRADS is its successful solving of a matrix completion problem with $15,694,167$ constraints and a matrix variable of size $40,000 \times 40,000$ in $351$ seconds. 

\end{abstract}

\noindent {\bf Keywords:}
semidefinite programming, alternating direction method of multipliers, low-rank factorization, logarithmic rank reduction, large-scale SDP
\section{Introduction}
In this paper, our goal is to address the general semidefinite program (SDP)  in the following form:
\begin{equation}
\min_{X\in\mbb S^{n}} \ \ f(X)\ \ 
\st  \ \ \mcal A(X)=b, \ \ X\succeq 0,
\label{prob:SDP}
\end{equation}
where $\mathbb{S}^{n}$ represents the space of all the  $n \times n$ symmetric matrices, $f(X)$ denotes a convex and continuously differentiable function, and the linear map $\mathcal{A}(\cdot):\mathbb{S}^n\rightarrow\mbb {R}^m$ is defined as $\mathcal{A}(X)=\left(\langle A_1, X\rangle, \ldots, \langle A_m, X\rangle\right)^{\top}$ with each $A_i\in \mathbb{S}^n$ for $i=1,2,\ldots, m$. Specifically, when $f(X)$ has the linear form, we have the linear SDP problem:
\begin{equation}
\min_{X\in\mbb S^{n}} \ \ \inner{C}{X}\ \ 
\st  \ \ \mcal A(X)=b, \ \ X\succeq 0,
\label{prob:linear-SDP}
\end{equation}
where $C\in \mbb S^n$. Throughout this paper, we conduct theoretical analyses on the general SDP form \eqref{prob:SDP} while focusing on the linear form \eqref{prob:linear-SDP} for computational explorations.

SDPs have a wide range of applications in practice, such as the MaxCut problem~\citep{goemans1995improved}, optimal power flow~\citep{lavaei2011zero}, combinatorial optimization~\citep{boyd1997semidefinite} and sensor network localization~\citep{so2007theory}, to name a few.
Traditionally, interior point methods are the go-to solution for solving SDP problems. It has been well known that for SDP, interior point methods can achieve high accuracy in polynomial time.  See ~\cite{alizadeh1995interior, nesterov1994interior, alizadeh1991combinatorial}.
Despite their effectiveness, traditional interior point methods struggle with scalability for large-scale problems due to the necessity of solving Newton systems repeatedly, which is computationally intensive.
To bypass the Newton system and enhance scalability, more recent research has shifted towards the first-order methods~\citep{wen2010alternating, lan2011primal, odonoghue2016conic}.
A significant area of interest has been the use of penalty methods, such as the augmented Lagrangian method and alternating direction method of multipliers (ADMM), for solving SDPs~\citep{wen2010alternating,zhao2010newton,yang2015sdpnal+,wang2023decomposition}. Specifically, \cite{zhao2010newton}  introduced SDPNAL, which employs the augmented Lagrangian method on dual SDP and uses an inexact semi-smooth Newton-CG method to solve the inner subproblem.  
The subsequent work~\citep{yang2015sdpnal+} presented an enhanced algorithm SDPNAL+, which achieves improved efficiency and robustness and adequately addresses the numerical difficulty of SDPNAL when solving degenerate SDPs. 
\cite{odonoghue2016conic} proposes a general splitting conic solver (SCS) for linear and conic programming, which applies ADMM to solve the homogeneous self-dual embedding associated with SDP. 
For large-scale SDP, a critical issue is the quadratically growing ($\mathcal{O}(n^2)$) memory cost in maintaining the matrix variable.  
For a class of rank-constrained matrix problems, 
\cite{yurtsever2017sketchy} apply the conditional gradient method to their convex relaxation problems. To reduce the storage overheads, they employ random sketching to keep track of a low-rank approximation of the matrix variable.
In a follow-up work by \cite{yurtsever2021scalable},
random sketching is further incorporated in 
the augmented Lagrangian-based method to deal with more general trace-constrained SDPs.

In addition to the mentioned convex approaches, 
another important research direction for large-scale SDP is to exploit the low-rank structure of SDP solutions directly through nonconvex reformulation. 
In the work of \citet{pataki1998rank, barvinok1995problems}, the authors have shown that standard SDP admits an optimal solution with rank $r$ satisfying $r(r+1)/2\leq m$ for specific $r>0$.  
Due to this low-rank structure, \cite{burer2003nonlinear} propose to replace the semidefinite matrix $X$ with the factorized form $X=UU^{\top}$, where  $U\in {\mbb {R}}^{n\times r}$. This factorization gives rise to the following nonconvex formulation of \eqref{prob:linear-SDP}:
\begin{equation}
\min_{U\in\mbb {R}^{n\times r}} \ \ \inner{C}{UU^{\top}}\ \ 
\st \ \ \mcal A(UU^{\top})=b.
\label{prob:linear-BM}
\end{equation}
\cite{burer2003nonlinear, burer2005local} propose an effective limited memory BFGS augmented Lagrangian algorithm to handle the low-rank formulation.
Another popular low-rank approach for SDP is the Riemannian optimization methods~\citep{journee2010low, wen2013feasible, boumal2015riemannian}. These methods are built upon the assumption that the feasible set of \eqref{prob:linear-BM} is a Riemannian manifold.  Within the low-rank formulation, \cite{tang2023feasible}  proposes a rank-support adaptive Riemannian optimization method for solving the general convex SDP problems and gives the convergence rate analysis. There are also several different ways to solve the low-rank formulation of SDP. \cite{bellavia2021relaxed} propose a relaxed variant of the interior point method, which exploits the low-rank structure of the optimal solution. \cite{erdogdu2022convergence} apply the block-coordinate ascent method to solve the low-rank reformulation of SDP with diagonal constraints and establish global convergence towards the first-order stationary point. 
In the BM formulation, a critical step is determining the rank of the matrix variable. The rank $r$ is chosen to ensure that either an optimal solution for \eqref{prob:linear-BM} is attainable or a satisfactory approximate low-rank solution is obtained.
In particular, 
\cite{soandye2008} demonstrate the existence of a solution of logarithmic rank order that still maintains good approximation property. 
This discovery forms the basis of an empirical tactic to initiate the algorithm with a solution of rank $\mathcal{O}(\log(m))$, which will be employed in our approach. 

Our paper's contributions are summarized as follows: First, we introduce a novel ADMM algorithm tailored for general SDPs that accommodate nonlinear and convex objectives. By introducing a new matrix splitting approach, we recast the SDP into a bilinearly constrained and quadratically penalized nonconvex problem, much different from the Burer-Monteiro approach. This formulation offers two primary benefits for ADMM implementation: 1) For a linear objective,  enabled by the splitting technique, both the alternating steps in ADMM result in a convex 
 and well-conditioned quadratic problem, and hence, can be efficiently solved by either direct methods or indirect ones such as the conjugate gradient descent method. 
 This simplification provides a substantial advantage over the nonconvex quartic subproblems inherent in the BM formulation;
2) This splitting provides an enlarged search area, which could potentially improve the algorithm's performance on some hard instances.
Under mild conditions, we prove that the ADMM algorithm achieves a sublinear rate of convergence to the KKT solutions.  
Additionally, when the \L{}ojasiewicz inequality is satisfied with an exponent of ${1}/{2}$,  we enhance the convergence result to a linear rate, offering a significant improvement in the algorithm's efficiency.

Second, we present a new SDP solver, which employs low-rank factorization and ADMM through matrix splitting (\oursolver{}). This solver integrates the proposed ADMM along with several new implementation strategies that significantly improve the practical efficiency. 
Our empirical evaluations highlight two critical heuristics: 1) a warm-start procedure employing the Burer-Monteiro algorithm, and 2) aggressive rank reduction to logarithmic order with a dynamic rank selection strategy ensuring the stability.
By leveraging these enhancements, we compare \oursolver{} with the leading SDP solvers on a wide range of problems, including MaxCut, matrix completion, and general SDP benchmarks in SDPLIB and Mittlemann's dataset. Extensive numerical studies demonstrate that \oursolver{} is highly competitive and often outperforms the popular SDP solvers. Furthermore, \oursolver{} demonstrates notable efficacy in addressing super large-scale problems, an area where current SDP solvers show limitations. A case in point is \oursolver{}'s ability to solve a matrix completion problem with $15,694,167$ constraints and a $40000 \times 40000$ matrix variable within approximately 350 seconds.

A closely related work is \cite{chen2023burer},  which also propose an ADMM-based method modified from the {\BM} factorization. 
Similar to our work, they employ matrix splitting to generate more manageable subproblems. However, their focus is restricted to a particular subset of SDP problems where the linear constraint is only applied to the diagonal elements. Due to the diagonally-constrained structure, when converted to the BM formulation, the closed-form solution of the projection can be efficiently computed.
In stark contrast, our study encompasses a much wider range of SDP problem types, incorporating general linear constraints alongside convex nonlinear objectives. 
Given the challenges associated with handling general SDP, we adopt a substantially different approach in our convergence analysis to achieve convergence and complexity results.

The rest of our paper is organized as follows. Section~\ref{sec:admm} introduces the proposed ADMM algorithm with matrix splitting for semidefinite programming. Section~\ref{sec:convergence} establishes the convergence and convergence rate of the proposed algorithm. The missing proof can be found in Appendix~\ref{sec:theory}. Section~\ref{sec:solver} develops an enhanced solving approach based on ADMM and several useful implementation heuristics. 
Section~\ref{sec:numerical} examines the empirical performance of our solver on various SDP problems. Finally, Section~\ref{sec:conclusion} concludes our paper and points out our future direction.

\paragraph{Notations.} Let $\|\cdot \|_2$ denote the Euclidean norm for vectors and operator norm for matrices, and $\norm{\cdot}_F$ denote the Frobenius norm for matrices. For a given matrix $X$ and radius $\delta$, we use $\mcal B_{\delta}\brbra{X}$ to denote the ball of radius $\delta$ around $X$, i.e. $\mcal B_{\delta}\brbra{X}=\{Y|\norm{X-Y}_F\leq\delta\}$. For a nonempty and compact set $\mcal S$, we reuse the notation $\mcal B_{\delta}\brbra{\mcal S}$ to denote the neighborhood of $\mcal S$, which is defined as $\mcal B_{\delta}\brbra{\mcal S}=\{Y|\inf_{Z\in \mcal S} \norm{Y-Z}_F\leq \delta\}$. For any matrix $X$, let $\sigma_{\min}(X)$ denote the minimal singular value of $X$. Let $\norm{\mcal A}_2$ denote the operator norm of the linear map $\mathcal{A}$ and $s_{A}=\sum_{i=1}^{m}\norm{A_{i}}_F^{2}$. Given a sequence of matrices $\{X^k\}_{k=0}^{K-1}$ and a sequence of vectors $\{x^k\}_{k=0}^{K-1}$, we use $\Delta X_k$ to denote $X^k-X^{k-1}$ and $\Delta x^k$ to denote $x^k-x^{k-1}$, respectively.

\section{The Low-Rank Bilinear Decomposition and the ADMM Approach\label{sec:admm}}
In this section, we present a general low-rank ADMM method with matrix splitting for solving the SDP problem. While we shall develop a more practical SDP solver for the standard linear SDP~\eqref{prob:linear-SDP} in Section \ref{sec:solver}, the ADMM algorithm and the subsequent analysis are adaptable to accommodate a general convex objective beyond the linear function. For this reason, this and the next sections will focus on the general convex SDP \eqref{prob:SDP}.

\subsection{Bilinear Reformulations and The Penalized Problem}
The {\BM} factorization for the problem is expressed as follows:
\begin{equation}
\min_{U\in \mbb {R}^{n\times r}} \ \ f({UU^{\top}})\ \ 
\st \ \ \mcal A(UU^{\top})=b.
\label{prob:BM}
\end{equation}
To address the challenges posed by large-scale SDP instances and to thoroughly explore the low-rank domain, we propose a different and bilinear factorization technique for problem \eqref{prob:SDP}, yielding the following formulation:
\begin{equation}
\min_{U,V\in\mbb {R}^{n\times r}}  \ \ f(UV^{\top})\ \ 
\st \ \ \mcal A(UV^{\top})=b,\ \ 
\ \ U=V.
\label{prob:split}
\end{equation}
\noindent We further consider the penalized version of \eqref{prob:split} by treating $U-V$ as a penalty term in the objective 
\begin{equation}
\min_{U,V\in\mbb {R}^{n\times r}}  \ \ f(UV^{\top})+\frac{\gamma}{2}\norm{U-V}_F^2\ \ 
\st \ \ \mcal A(UV^{\top})=b.
\label{prob:penalty}
\end{equation}
\noindent 
To elucidate the relationship between model \eqref{prob:BM} and \eqref{prob:penalty} with a suitably large parameter $\gamma$, we introduce the following assumption.
\begin{assumption}
\label{assu:Lipschitz}
The set collected by the KKT points of \eqref{prob:penalty}, denoted by $\mcal W$, is nonempty.
The function $f$ is symmetric in $\mbb {R}^{n\times n}$, i.e., $f(X)= f(X^\top)$ for any matrix $X\in \mbb {R}^{n\times n}$, and $\nabla f$ and $\nabla^2 f$ are Lipschitz continuous on any convex compact set.
\end{assumption}
For any given $\sigma>0$ and $X\in \mcal B_{\sigma^2}\brbra{0^{n\times n}}$,  there exist constants $L_f^\sigma$ and $L^\sigma$ such that
$\|\nabla f(X)\|_F\leq L_f^\sigma$ and
 $\nabla f$ is $L^\sigma$-Lipschitz continuous. Moreover, $\nabla^2 f$ is $L_0^\sigma$-Lipschitz continuous on 
$\mcal B_{\sigma^2}\brbra{0^{n\times n}}$, i.e. 
\[
\|\nabla^2 f(X)[H]-\nabla^2 f(Y)[H]\|_F\leq L_0^\sigma
\|H\|_F \|X-Y\|_F,\quad \forall X,Y\in \mcal B_{\sigma^2}\brbra{0^{n\times n}}, H\in \mbb {R}^{n\times n},
\]
where $\nabla^2 f(X)[H]$ denotes the tensor product of $\nabla^2 f(X)$ and $H$.
From the symmetry of $f$, it is easy to verify that the partial gradients satisfy
\[
\frac{\partial {f(UV^\top)}}{\partial U}=\nabla f(UV^\top) V\quad{\rm and}\quad
\frac{\partial {f(UV^\top)}}{\partial V}=\nabla f(VU^\top) U.
\]
With this preparation, we are ready to give the definitions of $\epsilon$-KKT point of problems \eqref{prob:BM} and \eqref{prob:penalty}.
\begin{definition}
    We say that $(U^*,\lambda^*)$ is an $\epsilon$-KKT point of problem \eqref{prob:BM} if
    \begin{equation}
        \begin{split}
        \Bnorm{2\nabla f(U^*{(U^*)}^{\top})U^* + \sum_{i=1}^m \lambda_i^* A_i U^*}_F\leq \epsilon,\\
        \bnorm{\mcal A(U^*{(U^*)}^{\top})-{b}}_2\leq\epsilon,
        \end{split}
        \label{def:eKKT}
    \end{equation}
and $(\overline{U},\overline{V},\overline{\lambda})$ is an $\epsilon$-KKT point of problem \eqref{prob:penalty} if it satisfies
    \begin{equation}
     \begin{split} \Bnorm{\nabla f(\overline{U}\overline{V}^{\top})\overline{V} + \gamma(\overline{U}-\overline{V}) + \sum_{i=1}^m \overline{\lambda}_i A_i \overline{V}}_F\leq \epsilon,\\
        \Bnorm{\nabla f(\overline{V}\overline{U}^{\top})\overline{U} + \gamma(\overline{V}-\overline{U}) + \sum_{i=1}^m \overline{\lambda}_i A_i \overline{U}}_F\leq \epsilon,\\
        \bnorm{\mcal A(\overline{U}\overline{V}^{\top})-{b}}_2\leq\epsilon.
        \end{split}
         \label{def:eKKT-split}
    \end{equation}
    When $\epsilon=0$, the definition of $\epsilon$-KKT points reduces to that of KKT points.
\end{definition}

The following theorem states that the solution produced by our algorithm for~\eqref{prob:penalty} can be suitably adapted to address~\eqref{prob:BM} approximately.

\begin{theorem}
\label{prop:KKT}
    Let $(\overline{U}, \overline{V}, \overline{{\lambda}})$ be any $\epsilon$-KKT point of problem \eqref{prob:penalty} in
     $\mcal B_{\sigma}\brbra{0^{n\times r},0^{n\times r},0^{m}}$
    and $\widehat{U} = \frac{\overline{U} + \overline{V}}{2}$. Suppose Assumption~\ref{assu:Lipschitz} holds, then the following statement s valid.
    \begin{equation}
    \begin{split}
        & \left\| {\cal A}(\widehat{U}{\widehat{U}}^{\top}) - b \right\|_2 \leq \epsilon + \frac{\|\mathcal{A}\|_2(\epsilon+L_f^\sigma \sigma+\sqrt{s_A}\sigma^2)^2
}{4\gamma^2}, \\
        & \left\| 2\nabla{f(\widehat{U}\widehat{U}^{\top}) }\widehat{U}+2\sum_{i=1}^{m}\overline{\lambda}_iA_{i}\widehat{U} \right\|_F \leq 2 \epsilon+ \frac{(2L_0^\sigma \sigma^3 +3L^\sigma \sigma)(\epsilon+L_f^\sigma \sigma+\sqrt{s_A}\sigma^2)^2}{ 2\gamma^2}. 
    \end{split}
    \end{equation}
   Moreover, suppose that we choose the penalty parameter $\gamma$  such that 
\[
\gamma >\max\left\{  \frac{\sqrt{\|\mathcal{A}\|_2}(\epsilon+L_f^\sigma \sigma+\sqrt{s_A}\sigma^2)
}{2\sqrt{2\epsilon}},
 { \frac{\sqrt{2L_0^\sigma \sigma^3 +3L^\sigma \sigma}(\epsilon+L_f^\sigma \sigma+\sqrt{s_A}\sigma^2) }{ \sqrt{2 \epsilon}}}\right\},
\]
then $(\widehat{U},2\overline{\lambda})$ is a
$3\epsilon$-KKT point of problem~\eqref{prob:BM}.
\end{theorem}
Theorem~\ref{prop:KKT} further indicates that if we have an $\epsilon$-KKT point of the penalty problem with $\gamma=\Omega(1/\sqrt{\epsilon})$, then we are able to obtain an $\epsilon$-KKT point of problem~\eqref{prob:BM}. As we will see later, this theoretical insight corroborates the empirical evidence observed in our study.

\subsection{The ADMM Approach}
We describe ADMM for solving the penalty problem~\eqref{prob:penalty} in Algorithm \ref{alg:ADMM-for-SDP}.
In this algorithm, the  augmented Lagrangian function of \eqref{prob:penalty} is defined by
\begin{equation}
\mcal L_{\rho}(U,V,{\lambda})=f(UV^{\top})+\frac{\gamma}{2}\norm{U-V}_{F}^{2}+\inprod{{\lambda}}{\mcal A(UV^{\top})-b}+\frac{\rho}{2}\norm{\mcal A(UV^{\top})-b}_2^{2}.\label{eq:ALF}
\end{equation}
By leveraging the bilinear factorization, ADMM alternates between optimizing two proximal subproblems associated with $U$ and $V$.
Since $f$ is convex, both the subproblems are $\gamma$-strongly convex with their own variables. The optimality conditions of the subproblems can be written as
\begin{equation}
\begin{split}0 = & \nabla f(U^{k+1}(V^{k})^{\top})V^{k}+\gamma(U^{k+1}-V^{k})+\sum_{i=1}^{m}\lambda_{i}^{k}A_{i}V^{k}\\
&+\rho\sum_{i=1}^{m}\left(\inprod{A_{i}V^{k}}{U^{k+1}}-b_{i}\right)A_{i}V^{k},
\end{split}
\label{eq:optU}
\end{equation}
and
\begin{equation}
\begin{split}0 = & \nabla f(V^{k+1}(U^{k+1})^{\top})U^{k+1}-\gamma(U^{k+1}-V^{k+1})+\sum_{i=1}^{m}\lambda_{i}^{k}A_{i}U^{k+1}\\
&+\rho\sum_{i=1}^{m}\big(\inprod{A_{i}U^{k+1}}{V^{k+1}}-b_{i}\big)A_{i}U^{k+1}.
\end{split}
\label{eq:optV}
\end{equation}
For linear SDPs, it is noteworthy to point out that due to the $\gamma$ penalty term from the splitting, both the subproblems (\ref{eq:optU}) and (\ref{eq:optV}) reduce to well-conditioned linear systems, which can be solved effectively by the conjugate gradient (CG) method. Numerical results supporting this view can be found in Section~\ref{sec:cg}.

\ifthenelse{\boolean{Arxiv}}{
\begin{algorithm}[H]
\SetAlgoNoLine
\caption{\admm{}\label{alg:ADMM-for-SDP}}
\label{alg:abip-with-restart}
Initialize $U^0,V^0,\lambda^0$ and penalty parameter $\rho$. \\
Set $ k = 0$ \\
\While{the termination criteria is not satisfied}{
    $U^{k+1}\leftarrow \argmin_{U}\mcal L_{\rho}(U, V^k, \lambda^k)$ \\
    $V^{k+1}\leftarrow \argmin_{V}\mcal L_{\rho}(U^{k+1},V,\lambda^{k})$ \\
    $\lambda^{k+1} = \lambda^{k} + \rho\brbra{\mcal A(U^{k+1}(V^{k+1})^\top ) - b}$\\ 
    Set $k \leftarrow k + 1$
}
\end{algorithm}
}{}

\noindent For the stopping criterion of ADMM, we terminate the algorithm if
\begin{equation}\label{eq:stop}
\max \left\{\frac{\|\mathcal{A}(U^{k+1}(V^{k+1})^\top)-b\|_2}{1+\|b\|_\infty},\frac{\gamma\|V^{k+1}-V^{k}\|_F}{1+\|V^k\|_F}\right\} < \epsilon
\end{equation}
is met, where $||\cdot||_\infty$ is the infinity norm of a vector. From \eqref{eq:optU} and \eqref{eq:optV}, we can easily derive that
\begin{equation}
\begin{split} &\nabla f(U^{k+1}(V^{k+1})^{\top})V^{k+1}+\gamma(U^{k+1}-V^{k+1})+\sum_{i=1}^{m}\lambda_{i}^{k+1}A_{i}V^{k+1}\\
={} &\nabla f(U^{k+1}(V^{k+1})^{\top})(V^{k+1}-V^{k})+\left(\nabla f(U^{k+1}(V^{k+1})^{\top})-\nabla f(U^{k+1}(V^{k})^{\top})\right)V^{k}\\
&+\gamma(V^{k}-V^{k+1})+\rho\sum_{i=1}^{m}\inprod{A_{i}(V^{k+1}-V^{k})}{U^{k+1}}A_{i}V^{k}+\sum_{i=1}^{m}\lambda_{i}^{k+1}A_{i}(V^{k+1}-V^k)
\end{split}
\label{eq:optU_revise}
\end{equation}
and
\begin{equation}
\begin{split}
\nabla f(V^{k+1}(U^{k+1})^{\top})V^{k+1}+\gamma(V^{k+1}-U^{k+1})+\sum_{i=1}^{m}\lambda_{i}^{k+1}A_{i}V^{k+1}=0.
\end{split}
\label{eq:optv_revise}
\end{equation}
Assume $\|V^k\|_F\leq \sigma$
and $(U^{k+1},V^{k+1},\lambda^{k+1})\in \left(\mcal B_{\sigma}\brbra{0^{n\times r}}\right)^2\times \mcal B_{\sigma}\brbra{0^{m}}$.
By taking the Frobenius norm on both sides of \eqref{eq:optU_revise} and using triangle inequality, it is not difficult to show
\begin{equation}
\begin{split} &\|\nabla f(U^{k+1}(V^{k+1})^{\top})V^{k+1}+\gamma(U^{k+1}-V^{k+1})+\sum_{i=1}^{m}\lambda_{i}^{k+1}A_{i}V^{k+1}\|_F\\
\leq & 
\left({L_f^\sigma + L^\sigma \sigma^2  +  \rho s_A\sigma^2+\sqrt{s_A}\sigma  }+\gamma\right)\|V^{k+1}-V^k\|_F\\
\leq & \left({\frac{L_f^\sigma + L^\sigma \sigma^2  +  \rho s_A\sigma^2+\sqrt{s_A}\sigma}{\gamma}}+1\right)(1+\sigma)\epsilon\triangleq c\epsilon.
\end{split}
\label{eq:u-kkt}
\end{equation}
This implies that the $k$-th iterate $(U^{k+1},V^{k+1},\lambda^{k+1})$
is a $(\max\left\{c,1+\|b\|_\infty\right\}\cdot \epsilon)$-KKT point of problem \eqref{prob:penalty}
and thus $\big(\frac{U^{k+1}+V^{k+1}}{2},2\lambda^{k+1}\big)$ is an $O(\epsilon)$-KKT point of  \eqref{prob:BM} by taking $\gamma=\Omega(1/\sqrt{\epsilon})$.

\subsection{The Power of ADMM}
\begin{figure}
    \centering
    \begin{minipage}[t]{0.44\columnwidth}
        \includegraphics[width = 1.1\textwidth]{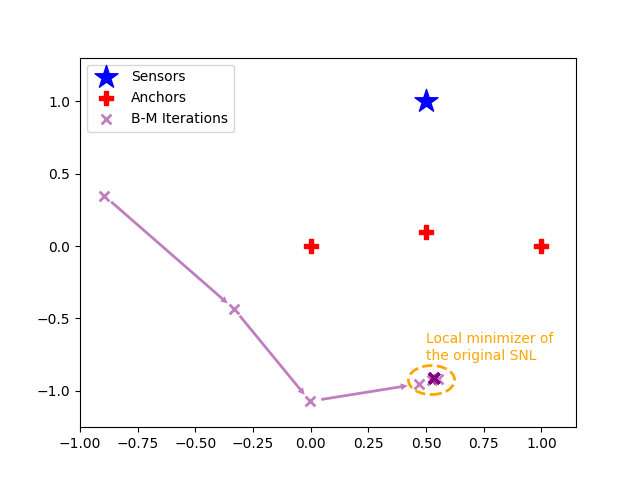}
    \end{minipage}
        \begin{minipage}[t]{0.44\columnwidth}
        \includegraphics[width = 1.1\textwidth]{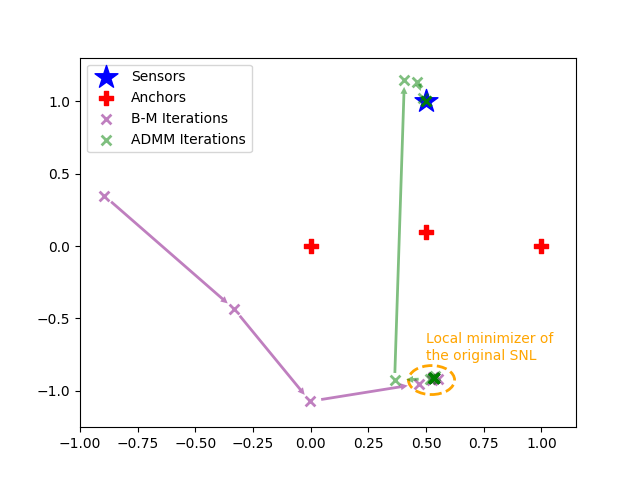}
    \end{minipage}
    \caption{\label{fig:escape}Comparison of ADMM and \BM{} Approaches on the SNL Instance}
\end{figure}

The intuition for constructing the bilinear factorization stems from the observation that within the iterative trajectory of \admm{}, $U$ and $V$ do not need to be identical. This splitting affords an enlarged search area, which sometimes improves the algorithm performance on hard instances. 
Figure~\ref{fig:escape} depicts a scenario in which the \BM{} approach is unable to solve the SDP relaxation of a sensor network localization (SNL) problem, whereas \admm{} demonstrates success. 
Given the location of anchors observed distances information between sensor-sensor and sensor-anchor pairs, the goal of SNL is to determine the sensor's true position, represented by the blue star. In a more formal setting, for a $d$-dimensional graph $G(V,E)$, we differentiate between the sensor set $V_x = \{x_1, x_2, \cdots, x_n\}$ and the anchor set $V_a = \{a_1, a_2, \cdots, a_m\}$ being a partition of $V$, where these sets collectively partition $V$ and $m \ge d + 1$. Observable edges are represented by two sets: $N_{x}= \{(i, j) \mid i, j \in V_x \}$ for sensor-sensor pairs and $N_{a}=\{(i, j) \mid i \in V_a, j \in V_x\}$ for sensor-anchor pairs. The lengths of these edges are denoted as $d_{i j} = \left\|x_{i}-x_{j}\right\|_2,~ \forall(i, j) \in N_{x}$ for sensor-sensor edges and $\tilde{d}_{k j} = \left\|a_{k}-x_{j}\right\|_2,~ \forall(k, j) \in N_{a}$ for sensor-anchor edges. The SNL problem is formulated as follows:

\begin{equation*}
\min_x \sum_{(i, j) \in N_{x}}\left|\left\|x_{i}-x_{j}\right\|_2^{2}-d_{i j}^{2}\right|+\sum_{(k, j) \in N_{a}}\left|\left\|a_{k}-x_{j}\right\|_2^{2}-\tilde {d}_{k j}^{2}\right|,
\end{equation*}
Let $e_{i} \in \mbb {R}^{n}$ be the vector with $1$ at $i$-th position and $0$ at others, and $\mathbf{0} \in \mbb {R}^{d}$ be the all-zeros vector. The SDP relaxation is then given by 
$$
\begin{array}{lll}
\max \ \ &0 \\
\text { s.t. } & Z_{(1: d, 1: d)}=I_{d} & \\
& \inprod {\left(\mathbf{0} ; e_{i}-e_{j}\right)\left(\mathbf{0} ; e_{i}-e_{j}\right)^{\top}}{Z}=d_{i j}^{2}, \quad \forall(i, j) \in N_x \\
& \inprod {\left(a_{k} ;-e_{j}\right)\left(a_{k} ;-e_{j}\right)^{\top}}{Z}=\tilde{d}_{k j}^{2}, \quad \forall(k, j) \in {N_a} \\
& Z \succeq 0,
\end{array}
$$
where the solution $Z \in \mbb {R}^{(d+n) \times(d+n)}$ comprises submatrices as:
$$
Z=\left[\begin{array}{cc}
I & X \\
X^{\top} & Y
\end{array}\right].
$$
If $Y = X^\top X$, then Z is a feasible solution for the above SDP problem, and $X = [x_1~ x_2~ \cdots~ x_n]$ gives a set of feasible locations of the sensors. 
The original SNL case under consideration involves a single sensor and three anchors, with their true positions illustrated in Figure~\ref{fig:escape}. We claim that the original SNL problem possesses a local minimum near the coordinates $(0.5, -0.9)$, posing potential difficulties for solving the SDP relaxation. The~\BM{} and~\admm{} methods' iterations, when mapped back to the original SNL context, indicate the sensor's coordinates in a two-dimensional space, as depicted in the figures. According to the left part of Figure~\ref{fig:escape}, the~\BM{} method becomes stuck at the local minimum of the original problem, thereby failing to place the sensor. For comparison, we initialize ADMM from the point where the~\BM{} method gets stuck. As shown on the right part of Figure~\ref{fig:escape}, ADMM successfully navigates away from this region and converges to the true position. The SNL instance is constructed as described in~\cite{lei2023blessing}, while the SDP relaxation of SNL  is detailed in~\cite{10.1145/984622.984630, 10.1145/1149283.1149286, doi:10.1137/060669395, shamsi2012sensor}.

\section{Convergence Analysis of \admm{}\label{sec:convergence}}

In this section, we present the main theoretical results of our proposed Algorithm~\ref{alg:ADMM-for-SDP}. We will prove the convergence of \admm{} for solving the bilinear factorized problem and further improve the convergence rate under the assumption of \L{}ojasiewicz inequality.

\subsection{Descent Property of \admm{}}

Our initial emphasis is on examining the change in the augmented Lagrangian function's value after each iteration. In the subsequent proposition, we demonstrate precisely how the updates of $U^{k+1}$, $V^{k+1}$, and ${\lambda}^{k+1}$ affect the augmented Lagrangian function's value.
\begin{proposition}\label{lem:decent} Suppose that Assumption~\ref{assu:Lipschitz} holds. Let $\bcbra{(U^k, V^k, \lambda^k)}_{k=0}^{K-1}$ be the sequence generated by Algorithm~\ref{alg:ADMM-for-SDP}, then for any $k\geq0$, we have 
\begin{align}
&\mcal L_{\rho}(U^{k+1},V^{k+1},{\lambda}^{k+1})-\mcal L_{\rho}(U^{k+1},V^{k+1},{\lambda}^{k})=\tfrac{1}{\rho}\norm{{\lambda}^{k+1}-{\lambda}^{k}}_2^{2}\label{eq:decrease1},\\
&\mcal L_{\rho}(U^{k+1},V^{k+1},{\lambda}^{k})-\mcal L_{\rho}(U^{k},V^{k},{\lambda}^{k}) \leq-\tfrac{\gamma}{2}(\norm{V^{k+1}-V^{k}}_F^{2}+\norm{U^{k+1}-U^{k}}_F^{2}).\label{eq:decrease2}
\end{align}
\end{proposition}

It is worth noting that the update of the dual variable ${\lambda}$ leads to an increase in the value of the augmented Lagrangian function. Therefore, bounding $\|\lambda^{k+1}-\lambda^{k}\|_2^2$ is crucial to ensure a sufficient decrease in each iteration. To this end, we need some boundedness and nonsingularity conditions.

\begin{assumption}
\label{assu:compact}
There exists a $\rho_{0}>0$ such that for any given $\beta$, the level set $S_{\beta}:=\{(U,V)\mid f(UV^{\top})+\frac{\rho_{0}}{2}\norm{U-V}_{F}^{2}+\frac{\rho_{0}}{2}\norm{{\cal A}(UV^{\top})-b}_2^{2}\leq\beta\}$
is either compact or empty.
\end{assumption}

We remark that Assumption~\ref{assu:compact}, while seemingly stringent, is applicable in a broad range of scenarios. The subsequent proposition demonstrates the widespread applicability of Assumption~\ref{assu:compact} across various contexts.

\begin{proposition}\label{propostion:assu1}
Assumption~\ref{assu:compact} holds for some $\rho_0>0$ if either of the following conditions is satisfied:
\begin{enumerate}
\item[(i)] $f(X)+\frac{\rho_0}{2}\|\mathcal{A}(X)-b\|_2^2$ strongly convex with respect to $X\in \mbb {R}^{n\times n}$.
\item[(ii)] $f(X)=\langle C, X\rangle$ with $C \in \mathcal{S}^n$, and the dual problem of~\eqref{prob:linear-SDP} has a strictly feasible solution $y^*$. Define $\lambda_{\min}$ and $\lambda_{\max}$ as the smallest and largest eigenvalues of $C-\sum_{i=1}^m y_i^*A_i$, respectively, and let $\kappa_c = \frac{\lambda_{\max}}{\lambda_{\min}}$ represent its condition number. Set $\rho_0 \geq 1 + \lambda_{\max} \kappa_c$.
\end{enumerate}
\end{proposition}

According to Assumption~\ref{assu:Lipschitz}, the nonlinear system $\mathcal{A}(UV^T) = b$ is feasible. Hence, we can locate an initial point $(U^0,V^0,\lambda^0)$ near the solution of $\mathcal{A}(UV^T)=b$ such that
\begin{equation}\label{initial}
\mcal L_{\rho}(U^{0},V^{0},{\lambda}^{0})+\frac{1}{\rho}\norm{{\lambda}^{0}}_2^{2}\leq \beta^0< \hat{\beta},
\end{equation}
where $\beta^0$ and $\hat{\beta}$ are constants independent of $\rho$. Using Assumption~\ref{propostion:assu1}, we establish that $S_{\hat{\beta}}$ is bounded, and there exists a constant $\delta > 0$ such that $\|U\|_F \leq \delta$ and $\|V\|_F \leq \delta$. 

We also need the following assumption, which assumes that the constraints $\mathcal{A}(UU^T) = b$ are linearly independent within the region $S_{\hat{\beta}}$. More previously, let the function $\mcal C:\mbb {R}^{n\times r}\rightarrow \mbb {R}^{nr\times m}$ denote the map from any point $U$ to the matrix consisting of the vectorized gradients of these constraints $\mathcal{A}(UU^T) = b$, i.e., $\mcal C(U)=[\vec{(}A_{1}U),\cdots,\vec{(}A_{m}U)]$. We assume that $\mcal C(U)$ is of full rank for any point $(U, V)$ in the region $S_{\hat{\beta}}$.
\begin{assumption}\label{assu:full rank}
There exists a constant $\eta>0$ such that 
\[
\min_{(U,V)\in S_{\hat{\beta}}}
\sigma_{\min}\brbra{\mcal C(U)}\geq\eta,
\]
where  and $\sigma_{\min}$ denotes the least singular value of a matrix.
\end{assumption}
Assumption~\ref{assu:full rank} is called the linear independence constraint qualification (LICQ), which is standard in the convergence analysis of ADMM-type methods. For example, see \cite{lin2015global,wen2010alternating,xie2021complexity}.

With these foundations in place, we are now prepared to establish the bound of $\|\lambda^{k+1}-\lambda^{k}\|_2^2$ and hence the descent property of our proposed Algorithm \ref{alg:ADMM-for-SDP}.

\begin{lemma}\label{lem:bound of dual variable} 
Suppose the initial point $(U^0,V^0,\lambda^0)$ is chosen according to \eqref{initial}, and  Assumptions~\ref{assu:Lipschitz}--\ref{assu:full rank} hold. Let $\rho$ and $\gamma$ be selected in Algorithm~\ref{alg:ADMM-for-SDP} such that 
\begin{equation}\label{ineq:rho:lowerbound}
\rho>\max\left\{2\rho_0,\frac{2c_{1}}{\gamma},\frac{2c_{2}}{\gamma},{\frac{1}{\hat{\beta}-\beta^0}} \left({\frac{L_{f}^{\delta}\delta+2\gamma \delta}{ \eta}}\right)^2\right\},\,\,\gamma >\rho_0.
\end{equation}
Then, for any $k\geq 0$, we have
\begin{itemize}
\item[(i)]\, $\mcal L_{\rho}(U^{k},V^{k},{\lambda}^{k})\leq \beta^0$ and $\mcal L_{\rho}(U^{k},V^{k},{\lambda}^{k})+\frac{3}{4\rho}\norm{{\lambda}^{k}}_2^{2}< \hat{\beta}$;
\item[(ii)]\,$(U^k,V^k)\in S_{\hat{\beta}}$;
\item [(iii)] $\norm{{\lambda}^{k+1}-{\lambda}^{k}}_2^{2}\leq c_{1}\norm{\Delta U^{k+1}}_F^{2}+c_{2}\norm{\Delta V^{k+1}}_F^{2},$
where $c_{1}=\frac{2}{\eta^{2}}\brbra{L^\delta\delta^{2}+L_{f}^\delta+\gamma+\frac{1}{\eta}L_{f}^\delta \delta\sqrt{s_{A}}+\frac{2}{\eta}\gamma \delta\sqrt{s_{A}}}^{2}$
and $c_{2}=\frac{2}{\eta^{2}}\brbra{L^\delta\delta^{2}+\gamma}^{2}$; 
\item [(iv)] $\mcal L_{\rho}(U^{k+1},V^{k+1},{\lambda}^{k+1})-\mcal L_{\rho}(U^{k},V^{k},{\lambda}^{k})\leq -(\frac{\gamma}{2}-\frac{c_{1}}{\rho})\norm{\Delta U^{k+1}}_F^{2}-(\frac{\gamma}{2}-\frac{c_{2}}{\rho})\norm{\Delta V^{k+1}}_F^{2}.$
\end{itemize}
\end{lemma}

\subsection{Convergence and Convergence Rate of ADMM}
Next, we proceed to establish the convergence of the whole sequence generated by Algorithm \ref{alg:ADMM-for-SDP}. To accomplish this, it is necessary to impose certain conditions on the augmented Lagrangian function, such as the \L{}ojasiewicz property. Denote the set collected by the KKT points of \eqref{prob:penalty} as $\mcal W$ and the description of \L{}ojasiewicz inequality is shown as follows.
\begin{assumption}[\L{}ojasiewicz Inequality] \label{assum:KL} 
    The augmented Lagrangian function $\mcal L_{\rho}(U,V,{\lambda})$ satisfies the \L{}ojasiewicz inequality with exponent $\alpha\in[0,1)$ at any KKT point of problem \eqref{prob:penalty}, i.e. there exist constants $c,\delta_{c}>0$ such that for any $(U,V,{\lambda})\in\mbb {R}^{n\times r} \times \mbb {R}^{n\times r} \times \mbb {R}^{m}$ and $(\overline{U},\overline{V},\overline{\lambda})\in \mcal W$ satisfying $\operatorname{dist}\brbra{(U,V,{\lambda}),(\overline{U},\overline{V},\overline{\lambda})}\leq\delta_{c}$,
    \begin{equation}\label{eq:KL inequality}
        \abs{\mathcal{L}_{\rho}(U,V,{\lambda})-\mathcal{L}_{\rho}(\overline{U},\overline{V},\overline{\lambda})}^{\alpha} \leq c \norm{\nabla \mathcal{L}_{\rho}(U,V,{\lambda})}_F.
    \end{equation}
\end{assumption}

Assumption~\ref{assum:KL} is non-restrictive as it remains valid in various instances. For example, when $f$ is a lower-semicontinuous semi-algebraic or tame function \citep{attouch2013convergence}, the augmented Lagrangian function $\mathcal{L}_{\rho}$ adheres to this condition. Furthermore, the \L{}ojasiewicz inequality inequality is satisfied for $\mathcal{L}_{\rho}$ at every interior point of its domain. With this assumption, we are ready to state the convergence of the sequence.

\begin{theorem}\label{thm:convergence}
    Suppose Assumptions~\ref{assu:Lipschitz}--~\ref{assum:KL} hold. Let the initial point $(U^0,V^0,{\lambda}^0)$ and $\rho,\gamma$ satisfy the conditions of Lemma~\ref{lem:bound of dual variable}. Then, the sequence of iterates $\{(U^k,V^k,{\lambda}^k)\}_{k\geq 0}^{K-1}$ converges to a KKT point of problem \eqref{prob:penalty}, $(\overline{U}, \overline{V},\overline{{\lambda}})$. Moreover, the total iteration complexity to achieve an $\epsilon$-KKT point is $\mcal O(\epsilon^{-1})$.
\end{theorem}

Therefore, we have established the $O(1/k)$ convergence rate of our ADMM method. However, as observed in our numerical experiments, Algorithm \ref{alg:ADMM-for-SDP} shows an even faster convergence rate in most cases. This is achievable if we could precisely characterize the exponent $\alpha$ of the \L{}ojasiewicz inequality.

\begin{theorem}\label{thm:linear convergence}
   Let the initial point $(U^0,V^0,{\lambda}^0)$ and $\rho,\gamma$ satisfy the conditions of Lemma \ref{lem:bound of dual variable}.
    Suppose that Assumption~\ref{assu:Lipschitz}--~\ref{assum:KL} hold. Then the following conditions hold
    \begin{enumerate}
        \item[(i)] If $\alpha=0$, then the sequence $\{(U^k,V^k,\lambda^k)\}_{k=0}^{\infty}$ converges in a finite number of steps.
        \item[(ii)] If $\alpha\in(0,\frac{1}{2}]$, then there exists $k_0\geq0$, $\xi>0$ and $\mu\in[0,1)$ such that for any $k>k_0$
        \begin{align*}
            \operatorname{dist}((U^k,V^k,\lambda^k),\mcal W)\leq \xi\mu^{k-k_0}.
        \end{align*}
        \item[(iii)] If $\alpha\in(\frac{1}{2},1)$, then there exists $k_0\geq0$, $\xi>0$ such that for any $k>k_0$
        \begin{align*}
            \operatorname{dist}((U^k,V^k,\lambda^k), \mcal W)\leq \xi (k-k_0)^{-\frac{1-\alpha}{2\alpha-1}}.
        \end{align*}
    \end{enumerate}
\end{theorem}

\section{An Enhanced Approach}\label{sec:solver}
In this section, we introduce \oursolver{}, an SDP solver based on the Low-rank ADMM Splitting approach, for effectively solving large-scale SDP problems. 
First, we give an outline of \oursolver{} in Section \ref{sec:loras}. Then, in Section~\ref{sec:instability}, we show the motivation of adding a warm-start phase that universally employing \admm{} does not guarantee stability. Subsequently, we propose a logarithmic dynamic rank selection strategy in Section~\ref{sec:4.3}, which significantly improves the performance of \oursolver{}. Finally, we detail the warm-started CG method used to solve the ADMM subproblems and show its good performance.

\subsection{\label{sec:loras}\oursolver{}}
\oursolver{}  is designed to derive a satisfactory low-rank approximate solution for the linear SDP problem~\eqref{prob:linear-SDP}. The approach unfolds in two phases: Phase I centers on the {\BM} factorization problem~\eqref{prob:BM}. To solve this problem, we employ the augmented Lagrangian algorithm~\citep{burer2003nonlinear}, initiating with a random starting point and employing the L-BFGS method to enhance algorithmic efficacy.

In the ensuing discussion, we refer to the augmented Lagrangian algorithm applied to solve problem~\eqref{prob:BM} as \alm{}.
When problem~\eqref{prob:BM} is solved to a certain accuracy, we stop the augmented Lagrangian algorithm and take the secured solution $(\tilde U^*, \tilde {\lambda}^*)$. Then, we switch to Phase II, focusing on solving the penalized bilinear decomposed problem (\ref{prob:penalty}) with \admm{}. We split the matrix variable $U$ and initiate \admm{} by setting $U^0 = V^0 = \tilde U^*, {\lambda}^0 = \frac{1}{2}\tilde {\lambda}^*$, where $(U^0, V^0, {\lambda}^0)$ is expected to be a satisfactory initial point as stated in Lemma \ref{lem:bound of dual variable}. The \admm{} is then executed to derive an approximate solution $(\overline{U}, \overline{V})$ for problem~(\ref{prob:penalty}) satisfying the final tolerance.
By setting $\widehat U = \frac{\overline{U}+\overline{V}}{2}$, we obtain an intermediate solution for problem (\ref{prob:BM}). Finally, with $\widehat X = \widehat U \widehat U^\top$, we get a solution to the general SDP problem (\ref{prob:SDP}). The complete pseudocode for the \oursolver{} is presented in Algorithm~\ref{alg:approach}.

\ifthenelse{\boolean{Arxiv}}{
\begin{algorithm}[H]
\SetAlgoNoLine
\caption{\oursolver{}\label{alg:approach}}
Initialize $\tilde U^0,\tilde {{\lambda}}^0$ and penalty parameter $\rho$\\
    \textbf{Phase I (the warm-start phase): }Call \alm{} to solve problem (\ref{prob:BM}) with initialization $(\tilde U^0,\tilde {{\lambda}}^0, \rho)$, and output a solution pair $(\tilde U^*,\tilde {{\lambda}}^*)$ once the \textbf{warm start} termination criteria is satisfied;\\ 
    \textbf{Switching Phases: }Set $U^0 = V^0 = \tilde U^*, {\lambda}^0 = \frac{1}{2} \tilde {{\lambda}}^*$;\\ 
    \textbf{Phase II (the solving phase): }Call \admm{} to solve problem (\ref{prob:penalty})  with initialization $(U^0,V^0,{{\lambda}}^0, \rho)$, and output a solution pair $(\overline U, \overline V)$ when the \textbf{final} termination criteria is satisfied;\\ 
    Set $\widehat U = \frac{\overline U + \overline V}{2}, \widehat X = \widehat U \widehat U^\top$ and output $\widehat X$ as a solution for (\ref{prob:SDP});
\end{algorithm}
}{}

\subsection{\label{sec:instability}The Motivation of the Two-phases Structure}
Although we aim at developing an approach based on the previously proposed \admm{}, our exploration reveals that directly applying Algorithm~\ref{alg:ADMM-for-SDP} to solve (\ref{prob:penalty}) might not always ensure stability.
Empirically, initiating from arbitrary starting points occasionally leads to suboptimal performance.  
In view of our theoretical result~(Theorem~\ref{thm:convergence}),  Algorithm~\ref{alg:ADMM-for-SDP} exhibits promising convergence properties when the initial point $(U^0, V^0, {\lambda}^0)$ aligns with the conditions specified in Lemma~\ref{lem:bound of dual variable}. 
These conditions imply that $U^0$ and $V^0$ should be close to the optimal solutions $U^*$ and $V^*$ to ensure fast convergence. 
To address this, we propose a warm-start approach, which initializes ADMM by the {\BM} method~\citep{burer2003nonlinear}, as BM keeps the low-rank factorization and directly yields a factorized solution $X = UU^\top$. To illustrate the ideas discussed above, we compare the warm-start approach against applying \admm{} throughout and applying \alm{} throughout on a set of representative problems selected\footnote{They are selected from the problems tested in Section \ref{sec:numerical}.}. In this experiment, we maintain a constant rank of $\sqrt{2m}$ and  employ SDPLR as the solver for \alm{}. Other implementation details and experiment settings follow Section \ref{sec:5.1} and \ref{sec:5.2}.

\ifthenelse{\boolean{Arxiv}}{
\begin{table}[!ht]
    \centering
    \caption{The effectiveness of The Warm-start Approach (Solving Times) }
    \begin{threeparttable}
    \begin{tabular}{c|ccc}
    \hline
        Problem & warm-started \admm{} & \admm{} & \alm{} (SDPLR) \\ \hline
        G60 & 3.47 & 6.61  & 7.66  \\ 
        MC\_3000 & 154 & t & 819  \\ 
        G40\_mb & 29.76 & 32.62  & 59.95  \\ 
        H3O & 19.37 & f & 199  \\ 
        p\_auss2\_3.0 & 56.85 & 158  & 83.95  \\ 
        qap7 & 1.41 & 2.14  & 2.18  \\ 
        qap10 & 6.47 & 21.51  & 10.19  \\ 
        theta12 & 62.26 & 22.76  & 127 \\ \hline
    \end{tabular}
    \begin{tablenotes}
\footnotesize
\item [1] We denote "timeout" by "t" and "fail" by "f", and the solving times are capped at 1800 seconds.
\item [2] For LRADMM, we use an initial $\rho = 100$.
\end{tablenotes}
\end{threeparttable}
\label{tab:warm-start}
\end{table}
}{}

From Table~\ref{tab:warm-start}, we observe that while \admm{} outperforms the \alm{} in solving problems such as G60, G40\_mb and theta12, it can be notably slower and may even fail in the other cases. However, the warm-started \admm{} exhibits a stable performance across the problems tested. 
Therefore, to enhance ADMM's stability, we integrate a {\WS} phase using \alm{}. 
More empirical results presented in Section~\ref{sec:numerical} indicate that using \alm{} for the {\WS} phase is adequately effective, though exploring more effective warm-start strategies presents an intriguing direction for future research.

\subsection{The Dynamic Rank Selection Strategy} \label{sec:4.3}
To improve the computation and memory efficiency, it is often desirable to apply a smaller rank while maintaining the algorithm's stability. Most literature on low-rank SDP algorithms typically recommends a rank of $\sqrt{2m}$ to maintain favorable properties of local minima. However, a theoretical result
by~\cite{soandye2008}  
suggests that an approximate solution of satisfactory quality can be efficiently achieved by employing a rank of $O(\log m)$.
The theorem underscoring this idea is stated as follows:
\begin{quote}
\textbf{Theorem 1.1 of~\citep{soandye2008}}
Let $A_{1}, \ldots, A_{m} \in \mbb {R}^{n \times n}$ be symmetric positive semidefinite matrices, and let $b_{1}, \ldots, b_{m} \geq 0$. Suppose that there exists an $X \succeq 0$ such that $\inprod{A_{i}}{X}=b_{i}$ for $i=1,2, \ldots, m$. Let $r=\min \{\sqrt{2 m}, n\}$. Then, for any $d \geq 1$, there exists an $X_{0} \succeq 0$ with $\operatorname{rank}\left(X_{0}\right) \leq d$ such that:
$$
\beta(m, n, d) \cdot b_{i} \leq \inprod{A_{i}}{X_{0}} \leq \alpha(m, n, d) \cdot b_{i} \quad \text { for } i=1, \ldots, m
$$
where:
$$
\alpha(m, n, d)= \begin{cases}1+\frac{12 \ln (4 m r)}{d} & \text { for } 1 \leq d \leq 12 \ln (4 m r)  \\ 1+\sqrt{\frac{12 \ln (4 m r)}{d}} & \text { for } d>12 \ln (4 m r)\end{cases}
$$
and
$$
\beta(m, n, d)= 
\begin{cases}\frac{1}{e(2 m)^{2 / d}} & \text { for } 1 \leq d \leq 4 \ln (2 m)  \\ 
\max \left\{\frac{1}{e(2 m)^{2 / d}}, 1-\sqrt{\frac{4 \ln (2 m)}{d}}\right\} & \text { for } d>4 \ln (2 m)\end{cases}
$$
Moreover, there exists an efficient randomized algorithm for finding such an $X_{0}$.
\end{quote}

This theoretical insight motivates us to aggressively adopt an $O(\log(m))$ rank to further reduce storage requirements and computational costs. 
Nevertheless, employing a rank $r$ that does not satisfy $r \geq \sqrt{2m}$ may make the problem more challenging. 
In some scenarios, fixing a rank of \(O(\log(m))\) may result in slow convergence or even failure to find a solution. Our strategy is simple: start with a small rank order (of \(O(\log(m))\)), increasing it only when necessary. We observed that adjusting the rank solely during the warm-start phase suffices, as   \admm{} is empirically more capable of handling low-rank instances than \alm{}. To determine when to increase the rank, we employ a heuristic to identify when a problem becomes "difficult" to solve, reflected by the number of inner iterations of the ALM subproblems. Upon detecting the difficulty, we halt the current subproblem's iterative cycle, increase the rank by a factor \(\alpha\), and then proceed with the algorithm's iterations. The results presented in Table~\ref{tab:rank} compare the performance of \oursolver{} when applying a consistent rank of \(r = 2\log(m)\), \(r = \sqrt{2m}\), or employing the dynamic rank selection strategy, where the rank is multiplied by a factor of 1.5 each time. We continue to use the selected problems in the previous subsection to illustrate the concept, with the experiment settings following the descriptions in Section~\ref{sec:5.1} and Section~\ref{sec:5.2}. From Table~\ref{tab:rank}, it can be seen that using a rank of 
$2\log(m)$ significantly reduces solving times compared to $\sqrt{2m}$ in many cases, 
while, for some other cases, a rank of $2\log(m)$ makes the problem slower to solve (e.g., qap7) or even unsolvable (e.g., qap10, theta12). On the other hand, the dynamic rank strategy can leverage the advantages of a low-rank structure while bypassing its potential challenges.

\ifthenelse{\boolean{Arxiv}}{
\begin{table}[!ht]
    \centering
        \caption{Solving Times Comparison of Different Rank Strategy}

    \begin{tabular}{ccc|cc|cc|cc}
    \hline
\multicolumn{3}{c|}{Problem} & \multicolumn{2}{c|}{$2\log(m)$} & \multicolumn{2}{c|}{$\sqrt{2m}$} & \multicolumn{2}{c}{Dynamic Rank} \\ \hline
        Name & $n$ & $m$ & Rank & Time & Rank & Time & Final Rank & Time \\ \hline
        G60 & 7000 & 7000 & 18 & 0.70  & 118 & 3.47  & 18 & 0.66 \\ 
        MC\_3000 & 3000 & 930328 & 27 & 9.94  & 1364 & 154  & 27 & 10.34 \\ 
        G40\_mb & 2001 & 2001 & 15 & 24.86  & 63 & 29.76  & 15 & 19.89 \\ 
        H3O & 3163 & 2964 & 16 & 10.14  & 77 & 19.37  & 16 & 10.12 \\ 
        p\_auss2\_3.0 & 9116 & 9115 & 18 & 4.45  & 135 & 56.85  & 18 & 4.42 \\ 
        qap7 & 50 & 358 & 12 & 9.03  & 27 & 1.41  & 26 & 1.02 \\ 
        qap10 & 101 & 1021 & 14 & f & 45 & 6.47  & 45 & 2.97 \\ 
        theta12 & 601 & 90020 & 23  & f & 424  & 62.26  & 77 & 86.23 \\ \hline
    \end{tabular}
        \label{tab:rank}
\end{table}
}{}

\subsection{Solving the ADMM Subproblems with the Warm-started CG} \label{sec:cg}
Since we consider the linear SDP problems, where the subproblems~\eqref{eq:optU} and~\eqref{eq:optV} of \admm{} with respect to $(U^{k+1}, V^{k+1})$ can be reformulated into the following linear systems with dimension $nr \times nr$:
\begin{equation}
\begin{split} &\left ( \rho \sum_{i=1}^{m}\mathrm{vec}(A_{i}V^{k})\mathrm{vec}(A_{i}V^{k})^\top+\gamma I_{nr\times nr} \right ) \mathrm{vec} (U^{k+1})\\ =& - \mathrm{vec} \left (CV^{k}-\gamma V^{k}+\sum_{i=1}^{m}\lambda_{i}^{k}A_{i}V^{k} - \rho \sum_{i=1}^{m} b_{i} A_{i}V^{k} \right )
,
\end{split}
\label{eq:linU}
\end{equation}
and
\begin{equation}
\begin{split} &\left ( \rho \sum_{i=1}^{m}\mathrm{vec}(A_{i}U^{k+1})\mathrm{vec}(A_{i}U^{k+1})^\top+\gamma I_{nr\times nr} \right ) \mathrm{vec} (V^{k+1})\\ =& - \mathrm{vec} \left (CU^{k+1}-\gamma U^{k+1}+\sum_{i=1}^{m}\lambda_{i}^{k}A_{i}U^{k+1} - \rho \sum_{i=1}^{m} b_{i} A_{i}U^{k+1} \right ),
\end{split}
\label{eq:linV}
\end{equation}
where $\mathrm{vec}(\cdot)$ is the vectorization of a matrix. 

To address the linear systems discussed, we employ a warm-started CG method. In this approach, we initialize the algorithm with $V^{k}$ (or $U^{k}$) to efficiently solve the linear system for $U^{k+1}$ (or $V^{k+1}$). It is important to highlight that the inclusion of the $\gamma$ penalty term from the splitting significantly improves the conditioning of the linear systems (\ref{eq:linU}) and (\ref{eq:linV}). Consequently, when $\gamma$ is sufficiently large, these systems are well-conditioned and can be effectively solved using the CG method without any preconditioning. To demonstrate the efficacy of this method, we present in Table \ref{tab:cgiter} the number of ADMM iterations, the total CG iterations required to solve the problem, and the average number of CG iterations per ADMM iteration (involving two subproblems) across the selection of representative problems. The problems are solved solely using \admm{} without \BM{} warm-start, and the ranks are set to be $\sqrt{2m}$. The result clearly shows that the average number of CG iterations per subproblem is extremely smaller than the worst-case scenario of $nr$ iterations.

\ifthenelse{\boolean{Arxiv}}{
\begin{table}[!ht]
    \centering
    \caption{Statistics of the Number of CG Iterations}
    \begin{threeparttable}

    \begin{tabular}{cc|ccc}
    \hline
        Problem & $nr =n \sqrt{2m}$ & ADMM Iterations & Total CG Iters & Avg CG Iters \\ \hline
        G60  & 828251 & 100 & 845 & 8.45  \\ 
        MC\_3000  & 4092176 & 17 & 9626 & 566.24  \\ 
        G40\_mb  & 126586 & 104 & 2252 & 21.65  \\ 
        p\_auss2\_3.0  & 1230829 & 1656 & 11152 & 6.73  \\ 
        qap7  & 1338 & 1829 & 91128 & 49.82  \\ 
        qap10  & 4564 & 29452 & 226521 & 7.69  \\ 
        theta12  & 255011 & 1261 & 6938 & 5.50 \\ \hline
    \end{tabular}
        \begin{tablenotes}
\footnotesize
\item [1] H3O is not included here since it fails to be solved using only \admm{} (see Table \ref{tab:warm-start}).
\end{tablenotes}

\end{threeparttable}
    \label{tab:cgiter}
\end{table}
}{}

\section{Numerical Experiments\label{sec:numerical}} In this section, we provide a detailed empirical evaluation of \oursolver{} applied to a diverse set of SDP problems. We begin by outlining the key implementation details and the experiment setup in Section~\ref{sec:5.1}. Subsequently, in Section~\ref{sec:5.2}, we assess the performance of our method on three sets of SDP problems: MaxCut, Matrix Completion, and a collection of general SDP benchmark problems.

\subsection{\label{sec:5.1}Implementation Details} 
We implement \oursolver{}\footnote{LoRADS is available at \url{https://github.com/COPT-Public/LoRADS}} in ANSI C with open-source numerical linear algebra libraries, Lapack~\citep{lapack99} and Blas~\citep{blackford2002updated}. To solve the subproblems of \admm{}, we use the warm-started conjugate gradient method to solve the linear systems \eqref{eq:linU} and \eqref{eq:linV}, where $f(UV^\top) = \inprod{C}{UV^\top}$. For the dynamic rank adjustment strategy discussed previously, we initiate the rank at $2\log(m)$ and increment the rank by a factor of 1.5 with each rank update. Especially, given the extensive scale of matrix completion problems from MC\_10000 to MC\_40000, a modified approach is adopted where the initial rank is set at $\log(m)$. Regarding the timing for variable splitting, we switch to Phase II upon meeting a switching tolerance of 1e-3 in the majority of instances tested. However, for all MaxCut problems (all problems in Section \ref{sec:5.3} and problems identified by starting with 'mcp' in Section \ref{sec:5.5}),
we apply a more lenient switching tolerance of 1e-2. Our code is capable of automatically selecting the phase-switching tolerance by analyzing the structural characteristics of the problem. The error used as the indicator for terminating Phase I is primal infeasibility, consistent with the final termination indicator outlined in Section \ref{sec:5.2}. Beyond these specific settings, our investigation suggests that the algorithm's effectiveness relies on applying certain techniques and implementation strategies. In the ensuing subsection, we delve into these aspects in detail.

\paragraph{The control of penalty factor.}
The efficacy of ADMM is significantly influenced by the choice of the penalty parameter, $\rho$. A well-known adaptive method for selecting this parameter is the residual balancing technique, introduced by \cite{he2000alternating, boyd2011}, whose idea is to balance primal and dual residuals. However, our experiments indicate that this method does not consistently provide satisfactory performance, primarily due to the challenge of selecting appropriate values for the related hyper-parameters across different problem instances. To avoid the need for problem-specific manual adjustment of the penalty parameter, we use a simple heuristic strategy that incrementally increases the penalty parameter within our approach. Specifically, within \admm{}, we increase $\rho$ by $1.2$ times every five iterations until it reaches a predetermined maximum value, $\rho_{\max}=5,000$. The initial value of  $\rho$ for \admm{} inherits from the warm start phase. In the warm-start phase, the update strategy for $\rho$ follows from  \cite{BMImplement}.

\paragraph{The heuristic factor.}
The appropriate value of $\rho$ for ADMM and  ALM may differ when transitioning between phases. To address this, we introduce a heuristic factor, $h$, to adjust $\rho$ appropriately during phase transitions. When switching from \alm{} to \admm{}, we multiply $\rho$ by  $h$ instead of carrying it over directly. In most cases, directly inheriting $\rho$ from the warm-start phase exhibits sufficiently robust performance. However, for certain types of problems, selecting an appropriate heuristic factor $h$ can significantly improve the performance. In subsequent experiments, we applied the heuristic factor as follows: $h = 10$ for all MaxCut, MaxCut MinBisection, graph partition, and quadratic assignment problems\footnote{The problems are recognized by: "G" and "mcp" for MaxCut; G40\_mb and G48\_mb for MaxCut MinBisection; "gpp" for graph partition; and "qap" for quadratic assignment problems.}; in matrix completion problems, we set $h=5$ for MC\_1000 to MC\_8000, and $h=2.5$ for MC\_10000 to MC\_40000; for all other problems, we set $h=1$ to indicate no use of the heuristic factor.

\paragraph{The implementation of \alm{}.}
To improve the performance of the \alm{} during the warm-start phase, we incorporate enhancements suggested in~\cite{BMImplement}. However, we do not apply the rank-updating method introduced in their work. Our implementation generally follows that of SDPLR v1.03\footnote{For the source code, see \url{https://yalmip.github.io/solver/sdplr/}. For a brief guide, see \url{https://sburer.github.io/files/SDPLR-1.03-beta-usrguide.pdf}.}.

\subsection{\label{sec:5.2}Experiment Setting}
We compare \oursolver{} with the following widely used SDP solvers:
\begin{itemize}
    \item SDPLR~\citep{BMImplement}. An open-source solver that utilizes the Burer-Monteiro (BM) low-rank factorization method.
    \item SDPNAL+~\citep
{sun2019sdpnal}. An open-source solver employing ALM, designed for tackling large-scale SDP problems.
    \item COPT~\citep
{ge2022cardinal}. A commercial solver based on the interior point method. COPT is currently recognized as the top-performing solver in Hans Mittelmann's benchmarks for SDP problems\footnote{\url{https://plato.asu.edu/ftp/sparse_sdp.html}.}.
This solver is included as a benchmark for the interior point method.
\end{itemize}

Given the diverse termination settings employed by different solvers—some of which are deeply integrated into their code—we have tried aligning their termination criteria while recognizing inherent differences. These solvers terminate when certain error metrics fall below a specified tolerance $\epsilon$:
\begin{itemize}
\item \oursolver{} and SDPLR terminate when $
\frac{\left\|\mathcal{A}\left(X\right) - b \right \|_2}{1+\left\|b\right \|_\infty} \leq \epsilon
$.

\item Denote the slack variable of the dual of \eqref{prob:SDP} by $S$ and the metric projection of $(X-S)$ onto $S_+^{n}$ by $\Pi_{S_+^{n}}(X-S)$, SDPNAL+ terminates when $
 \max \left\{
\frac{\|\mathcal{A}(X)-b\|}{1+\|b\|},\frac{\left\|\mathcal{A}^{*}(y)+S-C\right\|}{1+\|C\|},\frac{1}{5} \frac{\left\|X-\Pi_{S_+^{n}}(X-S)\right\|}{1+\|X\|+\|S\|}\right\}\leq \epsilon
$.

\item COPT terminates when $
\max \left\{\frac{\norm{\mcal A(X)-b}_2}{1+\norm{b}_{1}},\frac{\abs{\min\{0,\sigma_{\min}(C-\mcal A^{*}(\lambda))\}}}{1+\norm c_{1}}, \frac{\inner{C}{X}-\lambda^{\top}b}{1+\abs{\inner{C}{X}}+\abs{\lambda^{\top}b}} \right\} < \epsilon, $ where $c$ is the flatten vector of matrix $C$.

\end{itemize}

We have chosen to align the termination criterion of \oursolver{} with that of SDPLR rather than directly employing the condition~\eqref{eq:stop}. This decision is informed by the observation that the second term in \eqref{eq:stop}, $\frac{\gamma\|V^{k+1}-V^{k}\|_F}{1+\|V^k\|F}$, typically meets the required standards once the primal infeasibility gap is addressed. Notably, the primal infeasibility is measured as a value normalized using the infinity norm $\norm{b}_{\infty}$, making this approach significantly more rigorous than utilizing either $\norm{b}_2$ or $\norm{b}_{1}$, particularly in scenarios involving a large number of constraints.

In the following experiments in this section, we set $\epsilon$ to be 1e-5, and the time limit to be 10,000 seconds\footnote{The instances that failed to be solved in 10,000 seconds are marked with "t"}. The iteration cap for the LRADMM algorithm is set at 5,000. All the experiments are performed on a Macbook Pro with a 12 Core M3 Pro chip and 18 GB RAM. For clarity, this section only reports the results of solving times for each experiment. A report of properly defined and aligned solving errors for all solvers in all the experiments can be found in Appendix \ref{app:B}.

\subsection{The MaxCut Problems}\label{sec:5.3}

Given a weighted graph $G$ with weights $w$, the MaxCut problem aims to partition the vertices into two disjoint subsets such that the total weight of the edges connecting vertices across the subsets is maximized. 
The SDP formulation of the MaxCut problem can be expressed as follows:
\begin{align*}
    \max_{X \in {\mbb S^{n \times n}}} \ \ &\inprod{\frac{1}{4}L(G,w)}{X} \\ \st \ \ &X_{ii} = 1 \ \ , \forall i \in [n]\\
    &X \succeq 0,
\end{align*}
where the Laplacian matrix $L(G,w)$  is defined by
$$
L(G, w)_{i j}:= \begin{cases}-w_{i j} & \text { if } (i, j) \in E \\ \sum_{k} w_{i k} & \text { if } i=j \\ 0 & \text { o.w. }\end{cases}.
$$
We evaluated all previously discussed solvers using the Gset benchmark\footnote{For problem files, see \url{https://web.stanford.edu/~yyye/yyye/Gset/}}~\citep{ye_gset_2003}, and report the test results in Table~\ref{tab:max-cut-time}.   Our findings suggest that SDPNAL+ struggles with this specific set of SDP challenges, and 
\oursolver{} outperforms the other solvers in almost all cases.
For problems with $n,m < 5000$, \oursolver{} is approximately $1$ to $2$ times faster than both SDPLR and COPT. However, for larger instances where $n,m \ge 5000$, the performance gap widens significantly, with \oursolver{} being $10$ to $20$ times faster than SDPLR and COPT.
This trend highlights that the advantages of \oursolver{} become more pronounced in larger-scale problems. 

\ifthenelse{\boolean{Arxiv}}{
\begin{table}[!ht]
    \centering
    \small
    \caption{Solving Times on Max-Cut Problems}
    \begin{threeparttable}
    \setlength{\tabcolsep}{1.0mm}
    \begin{tabular}{c|c|cccc||c|c|cccc}
    \hline
        Problem & $n, m$ & \oursolver{} & SDPLR & SDPNAL+ & COPT & Problem & $ n, m$ & \oursolver{} & SDPLR & SDPNAL+ & COPT \\ \hline
        G1 & 800 & \textbf{0.18}  & 0.75  & 33.43  & 0.82  & G35 & 2000 & \textbf{0.23}  & 1.60  & 723  & 3.49  \\ 
        G2 & 800 & \textbf{0.17}  & 0.68  & 36.49  & 0.78  & G36 & 2000 & \textbf{0.24}  & 1.41  & 828  & 2.99  \\ 
        G3 & 800 & \textbf{0.18}  & 0.81  & 31.69  & 0.90  & G37 & 2000 & \textbf{0.25}  & 1.84  & 381  & 3.52  \\ 
        G4 & 800 & \textbf{0.17}  & 0.69  & 29.59  & 0.81  & G38 & 2000 & \textbf{0.25}  & 1.85  & 414  & 3.25  \\ 
        G5 & 800 & \textbf{0.15}  & 0.66  & 25.41  & 0.72  & G39 & 2000 & \textbf{0.38}  & 1.76  & 352  & 3.29  \\ 
        G6 & 800 & \textbf{0.23}  & 0.67  & 30.18  & 0.61  & G40 & 2000 & \textbf{0.57}  & 2.82  & 428  & 3.54  \\ 
        G7 & 800 & \textbf{0.22}  & 0.60  & 29.54  & 0.56  & G41 & 2000 & \textbf{0.52}  & 1.99  & 444  & 3.20  \\ 
        G8 & 800 & \textbf{0.20}  & 0.69  & 34.49  & 0.61  & G42 & 2000 & \textbf{0.54}  & 1.67  & 360  & 3.19  \\ 
        G9 & 800 & \textbf{0.20}  & 0.54  & 28.22  & 0.54  & G43 & 1000 & \textbf{0.09}  & 0.44  & 18.20  & 0.97  \\ 
        G10 & 800 & \textbf{0.20}  & 0.64  & 36.14  & 0.81  & G44 & 1000 & \textbf{0.11}  & 0.37  & 18.17  & 0.74  \\ 
        G11 & 800 & \textbf{0.05}  & 0.47  & 49.13  & 0.25  & G45 & 1000 & \textbf{0.10}  & 0.40  & 53.57  & 1.00  \\ 
        G12 & 800 & \textbf{0.05}  & 0.34  & 58.34  & 0.26  & G46 & 1000 & \textbf{0.10}  & 0.56  & 50.16  & 0.85  \\ 
        G13 & 800 & \textbf{0.05}  & 0.17  & 58.37  & 0.28  & G47 & 1000 & \textbf{0.10}  & 0.43  & 55.18  & 1.10  \\ 
        G14 & 800 & \textbf{0.07}  & 0.28  & 45.07  & 0.50  & G48 & 3000 & \textbf{0.12}  & 0.67  & 1436  & 2.09  \\ 
        G15 & 800 & \textbf{0.07}  & 0.34  & 55.99  & 0.59  & G49 & 3000 & \textbf{0.14}  & 0.76  & 462  & 2.00  \\ 
        G16 & 800 & \textbf{0.07}  & 0.33  & 56.73  & 0.58  & G50 & 3000 & \textbf{0.13}  & 0.96  & 1433  & 1.88  \\ 
        G17 & 800 & \textbf{0.06}  & 0.34  & 39.46  & 0.53  & G51 & 1000 & \textbf{0.09}  & 0.45  & 51.64  & 0.83  \\ 
        G18 & 800 & \textbf{0.12}  & 0.37  & 30.17  & 0.37  & G52 & 1000 & \textbf{0.09}  & 0.38  & 58.41  & 0.88  \\ 
        G19 & 800 & \textbf{0.12}  & 0.41  & 33.84  & 0.54  & G53 & 1000 & \textbf{0.09}  & 0.43  & 19.48  & 0.75  \\ 
        G20 & 800 & \textbf{0.12}  & 0.32  & 34.12  & 0.48  & G54 & 1000 & \textbf{0.08}  & 0.48  & 49.06  & 0.66  \\ 
        G21 & 800 & \textbf{0.11}  & 0.35  & 27.20  & 0.50  & G55 & 5000 & \textbf{0.38}  & 3.45  & 5587  & 18.04  \\ 
        G22 & 2000 & \textbf{0.28}  & 1.66  & 277  & 3.42  & G56 & 5000 & \textbf{0.42}  & 4.03  & 4664  & 18.18  \\ 
        G23 & 2000 & \textbf{0.31}  & 1.31  & 483  & 3.67  & G57 & 5000 & \textbf{0.54}  & 6.55  & 7123  & 11.34  \\ 
        G24 & 2000 & \textbf{0.27}  & 1.91  & 274  & 3.50  & G58 & 5000 & \textbf{0.87}  & 14.76  & 8426  & 24.52  \\ 
        G25 & 2000 & \textbf{0.27}  & 1.26  & 276  & 3.26  & G59 & 5000 & \textbf{2.42}  & 26.48  & 7318  & 25.16  \\ 
        G26 & 2000 & \textbf{0.27}  & 1.72  & 286  & 3.76  & G60 & 7000 & \textbf{0.66}  & 7.66  & t & 40.20  \\ 
        G27 & 2000 & \textbf{0.40}  & 1.66  & 525  & 3.49  & G61 & 7000 & \textbf{0.86}  & 6.69  & t & 41.16  \\ 
        G28 & 2000 & \textbf{0.50}  & 1.99  & 418  & 4.09  & G62 & 7000 & \textbf{1.02}  & 16.25  & t & 23.59  \\ 
        G29 & 2000 & \textbf{0.33}  & 2.37  & 272  & 4.27  & G63 & 7000 & \textbf{1.92}  & 23.69  & t & 59.32  \\ 
        G30 & 2000 & \textbf{0.31}  & 1.37  & 553  & 4.54  & G64 & 7000 & \textbf{5.60}  & 61.82  & t & 62.78  \\ 
        G31 & 2000 & \textbf{0.41}  & 2.03  & 321  & 3.86  & G65 & 8000 & \textbf{1.30}  & 24.80  & t & 33.35  \\ 
        G32 & 2000 & \textbf{0.18}  & 1.36  & 645  & 1.50  & G66 & 9000 & \textbf{1.44}  & 29.39  & t & 43.09  \\ 
        G33 & 2000 & \textbf{0.17}  & 1.21  & 475  & 1.44  & G67 & 10000 & \textbf{1.58}  & 42.73  & t & 57.18  \\ 
        G34 & 2000 & \textbf{0.17}  & 1.17  & 579  & 1.61 & ~ & ~ & ~\\ \hline
    \end{tabular}
\end{threeparttable}
\label{tab:max-cut-time}
\end{table}
}{}

\subsection{The Matrix Completion Problem}
The Matrix Completion Problem, which involves estimating the unobserved entries of a partially observed matrix, finds wide application in areas such as collaborative filtering, image inpainting, and system identification, among others. SDP is a well-known method for addressing matrix completion problems. However, the use of conventional SDP solvers in matrix completion, particularly in large-scale instances, is hindered by computational inefficiencies~\citep{li2019survey}. Despite these limitations, recent progress in SDP algorithms highlighted by~\cite{wang2023solving}, has demonstrated significant improvements over traditional SDP solvers. Additionally, recent studies by~\cite{yalccin2023semidefinite} have shown that the {\BM} factorization method is particularly advantageous in scenarios with limited observations and a small number of measurements. Motivated by these efforts, we test our method on matrix completion problems and demonstrate its effectiveness. 

Given a matrix $M \in \mbb {R}^{n \times m}$ with observed entries in $\Omega \subseteq \{1, \ldots, n \} \times \{ 1, \ldots, m \}$, the problem is to find $Y$ such that $Y_{ij} = M_{ij}$ for any $(i,j) \in \Omega$, while rank($Y$) being minimized. The convex optimization problem can be formulated as follows: 
\begin{align*}
\min_{Y \in \mbb {R} ^{m \times n}} \ \ \norm {Y}_*\ \ 
\st \ \  Y_{ij} = M_{ij}, \ \ \forall (i,j) \in \Omega,
\end{align*}
where $\norm {Y}_*$ is the nuclear norm of $Y$. Further, the nuclear norm minimization problem above can be formulated into a semidefinite programming problem
as follows:
\begin{align*}
\min_{Y \in \mbb {R} ^{m \times n}} \ \  \text{trace} \brbra{W_1} + \text{trace} \brbra{W_2} \ \ 
\st \ \  Y_{ij} = M_{ij}, \ \ \forall (i,j) \in \Omega, 
{\left[\begin{array}{cc}
W_1 & Y^{\top} \\
Y & W_2
\end{array}\right]} \succeq 0,
\end{align*}
which is equivalent to 
\begin{align*}
\min_{X \in \mbb S^{(m+n) \times (m+n)}} \ \   \inprod{I}{X},\ \ 
\st \ \  \bigg\langle\left[\begin{array}{cc}
0_{m \times m} & E_{i j}^{\top} \\
E_{i j} & 0_{n \times n}
\end{array}\right],{X}\bigg\rangle=2 M_{i j}, \ \ \forall (i,j) \in \Omega
, X \succeq 0,
\end{align*}
where $X = {\left[\begin{array}{cc}
W_1 & Y^{\top} \\
Y & W_2
\end{array}\right]}$, and $E_{ij}$ is the $ij$-th adjacency matrix with $(E_{ij})_{ij} = 1$ and 0 elsewhere.

\ifthenelse{\boolean{Arxiv}}{

\begin{table}[!ht]
    \centering
    \small
    \caption{Solving Times on Matrix Completion Problems}
    \begin{threeparttable}
    \begin{tabular}{c|cc|cccc}
    \hline
        Problem & $n$ & $m$ & \oursolver{} & SDPLR & SDPNAL+ & COPT \\ \hline
        MC\_1000 & 1000  & 199424  & \textbf{1.41}  & 36.96  & 6.21  & 14.84  \\ 
        MC\_2000 & 2000  & 550536  & \textbf{4.68}  & 220  & 60.65  & 106  \\ 
        MC\_3000 & 3000  & 930328  & \textbf{10.09}  & 819  & 343  & 416  \\ 
        MC\_4000 & 4000  & 1318563  & \textbf{18.78}  & 2401  & 326  & 1239  \\ 
        MC\_5000 & 5000  & 1711980  & \textbf{24.49}  & 4547  & 636  & 2502  \\ 
        MC\_6000 & 6000  & 2107303  & \textbf{32.15}  & 6390  & 2347  & 5006  \\ 
        MC\_8000 & 8000  & 2900179  & \textbf{57.21}  & t & 2921  & t \\ 
        MC\_10000 & 10000  & 3695929  & \textbf{27.33}  & t & 5879  & t \\ 
        MC\_12000 & 12000  & 4493420  & \textbf{49.40}  & - & - & - \\ 
        MC\_14000 & 14000  & 5291481  & \textbf{42.62}  & - & - & - \\ 
        MC\_16000 & 16000  & 6089963  & \textbf{74.66}  & - & - & - \\ 
        MC\_18000 & 18000  & 6889768  & \textbf{71.60}  & - & - & - \\ 
        MC\_20000 & 20000  & 7688309  & \textbf{83.83}  & - & - & - \\ 
        MC\_40000 & 40000  & 15684167  & \textbf{351}  & - & - & - \\ \hline
    \end{tabular}
        \begin{tablenotes}
\footnotesize
\item [1] For instances beyond MC\_10000, 
we only evaluate \oursolver{} and use "-" for other solvers as placeholders.
\end{tablenotes}
\end{threeparttable}
\label{tab:mc-time}
\end{table}
}

We benchmark our algorithm on a set of randomly generated matrix completion problems, adopting the generating process described by~\cite{wang2023solving}. 
The results, as depicted in Table~\ref{tab:mc-time}, demonstrate that \oursolver{} consistently outperforms other solvers across a range of instance sizes. For instances larger than MC\_10000, further evaluations of the compared solvers were not conducted due to time and computation resource limits. Notably, \oursolver{} successfully solves an exceptionally large-scale instance with $15,694,167$ constraints and a $40000 \times 40000$ matrix variable in less than 6 minutes.

\subsection{General SDP Benchmark Problems}\label{sec:5.5}
In this section, we evaluate our method across a wider range of SDP problems from SDPLIB\footnote{For the problem files and a detailed description, see \url{http://euler.nmt.edu/~brian/sdplib/}.} and Hans Mittelmann's SDP benchmarks\footnote{For the problem files and a detailed description, see \url{https://plato.asu.edu/ftp/sparse_sdp.html}.}. SDPLIB offers relatively small-scale SDP instances, while Hans Mittelmann's SDP benchmarks provide larger and more challenging SDP instances.

The comparative performance of the four solvers is detailed in Table~\ref{tab:sdplib-hans-time}. 
Our method exhibits considerable improvements over SDPLR across most tested instances.  In comparison with SDPNAL+, our solver falls short in only a few cases. Against COPT, our solver shows competitive performance, surpassing COPT in several instances, albeit falling short in others. We need to remark that the results of \oursolver{} are based on a very elementary parameter tuning detailed in Section \ref{sec:5.1}. For example, the problem "hand" can be solved in $5.81$ seconds if a more suitable $\rho_{max}$ and heuristic factor $h$ is applied, potentially reducing the solving time by a factor of 20.

\ifthenelse{\boolean{Arxiv}}{
\begin{table}[!ht]
\setlength{\tabcolsep}{3pt}
    \centering
    \footnotesize
    \caption{\label{tab:sdplib-hans-time}Solving Times on Problems from SDPLIB and Mittelmann's Set}
    \begin{threeparttable}
    \setlength{\tabcolsep}{0.7mm}
    \begin{tabular}{c|cc|cccc||c|cc|cccc}
    \hline
        Problem & $n$ & $m$ & \oursolver{} & SDPLR & SDPNAL+ & COPT & Problem & $n$ & $m$ & \oursolver{} & SDPLR & SDPNAL+ & COPT \\ \hline
        AlH$^\star$ & 5991 & 7230 & \textbf{63.21} & 2491  & 943  & 1208 & mcp250-1$^\divideontimes$ & 250 & 250 & \textbf{0.01} & \textbf{0.01} & 7.57 & 0.86 \\ 
        BH2$^\star$ & 2167 & 1743 & \textbf{7.80} & 124  & 48.31 & 35.65 & mcp250-2$^\divideontimes$ & 250 & 250 & \textbf{0.01} & \textbf{0.01} & 5.42 & 0.88 \\ 
        cancer\_100$^\star$ & 570 & 10470 & 513.88  & t & 728  & \textbf{83.36} & mcp250-3$^\divideontimes$ & 250 & 250 & \textbf{0.01} & 0.02 & 5.49 & 0.79 \\ 
        CH2$^\star$ & 2167 & 1743 & \textbf{6.74} & 51.46 & 88.84 & 19.06 & mcp250-4$^\divideontimes$ & 250 & 250 & \textbf{0.02} & 0.04 & 5.25 & 0.90 \\ 
        checker\_1.5$^\star$ & 3971 & 3971 & \textbf{0.95} & 40.45 & 6201  & 20.3 & mcp500-1$^\divideontimes$ & 500 & 500 & \textbf{0.01} & 0.03 & 21.36 & 0.20 \\ 
        G40\_mb$^\star$ & 2001 & 2001 & 19.89 & 59.95 & 4645  & \textbf{12.00} & mcp500-2$^\divideontimes$ & 500 & 500 & \textbf{0.02} & 0.04 & 14.59 & 0.20 \\ 
        G48\_mb$^\star$ & 3001 & 3001 & \textbf{17.21}  & 25.45 & 3019  & 17.50 & mcp500-3$^\divideontimes$ & 500 & 500 & \textbf{0.03} & 0.07 & 21.87 & 0.29 \\ 
        gpp100$^\divideontimes$ & 100 & 101 & \textbf{0.05} & 0.14 & 3.76 & 0.14 & mcp500-4$^\divideontimes$ & 500 & 500 & \textbf{0.04} & 0.11 & 23.22 & 0.40 \\ 
        gpp124-1$^\divideontimes$ & 124 & 125 & \textbf{0.04} & 0.26 & 6.14 & 0.24 & NH2$^\star$ & 2047 & 1743 & 24.77 & 19.03 & \textbf{10.81}  & 18.37 \\ 
        gpp124-2$^\divideontimes$ & 124 & 125 & \textbf{0.07} & 0.11 & 4.99 & 0.23 & NH3$^\star$ & 3163 & 2964 & \textbf{9.57} & 177  & 55.38  & 63.81 \\ 
        gpp124-3$^\divideontimes$ & 124 & 125 & \textbf{0.06} & 0.10 & 7.42 & 0.26 & p\_auss2\_3.0$^\star$ & 9116 & 9115 & \textbf{4.42} & 83.95 & t & 289 \\ 
        gpp124-4$^\divideontimes$ & 124 & 125 & \textbf{0.11} & 0.32 & 4.57 & 0.23 & qap5$^\divideontimes$ & 26 & 136 & 0.18 & 0.09 & 0.48 & \textbf{0.01} \\ 
        gpp250-1$^\divideontimes$ & 250 & 250 & \textbf{0.15} & 0.71 & 27.93 & 1.79 & qap6$^\divideontimes$ & 37 & 229 & 0.63 & 1.90 & 1.69 & \textbf{0.03} \\ 
        gpp250-2$^\divideontimes$ & 250 & 250 & \textbf{0.14} & 0.49 & 14.52 & 1.63 & qap7$^\divideontimes$ & 50 & 358 & 1.02 & 2.18 & 1.65 & \textbf{0.06} \\ 
        gpp250-3$^\divideontimes$ & 250 & 250 & \textbf{0.16} & 1.59 & 16.46 & 2.03 & qap8$^\divideontimes$ & 65 & 529 & 1.75 & 4.33 & 8.89 & \textbf{0.15} \\ 
        gpp250-4$^\divideontimes$ & 250 & 250 & \textbf{0.20} & 0.61 & 13.94 & 1.84 & qap9$^\divideontimes$ & 82 & 748 & 2.50 & 4.89 & 5.34 & \textbf{0.26} \\ 
        gpp500-1$^\divideontimes$ & 500 & 501 & \textbf{0.71} & 1.84 & 31.15 & 1.08 & qap10$^\divideontimes$ & 101 & 1021 & 2.97 & 10.19 & 6.14 & \textbf{0.43} \\ 
        gpp500-2$^\divideontimes$ & 500 & 501 & \textbf{0.72} & 1.02 & 106  & 0.76 & sensor\_500$^\star$ & 4029 & 3540 & 94.13 & 73.41 & 121  & \textbf{12.3} \\ 
        gpp500-3$^\divideontimes$ & 500 & 501 & \textbf{0.62} & 2.64 & 35.98 & 0.67 & shmup4$^\star$ & 800 & 4962 & \textbf{53.59} & 70.77 & 3068  & 115 \\ 
        gpp500-4$^\divideontimes$ & 500 & 501 & \textbf{0.65} & 1.30 & 36.37 & 0.76 & theta1$^\divideontimes$ & 50 & 104 & 0.22 & 0.32 & 0.32 & \textbf{0.02} \\ 
        H3O$^\star$ & 3163 & 2964 & \textbf{10.12} & 199  & 295  & 101 & theta2$^\divideontimes$ & 100 & 498 & 0.57 & 0.72 & 1.65 & \textbf{0.10} \\ 
        hand$^\star$ & 1297 & 1297 & 116  & \textbf{24.22} & 628  & 92.76 & theta3$^\divideontimes$ & 150 & 1106 & 1.49 & 2.29 & 1.57 & \textbf{0.30} \\ 
        ice\_2.0$^\star$ & 8114 & 8113 & \textbf{4.91} & 41.97 & t & 269 & theta4$^\divideontimes$ & 200 & 1949 & 3.09 & 6.01 & 3.74 & \textbf{0.82} \\ 
        mcp100$^\divideontimes$ & 100 & 100 & \textbf{$<$0.005} & \textbf{$<$0.005} & 1.19 & 0.07 & theta5$^\divideontimes$ & 250 & 3028 & 4.67 & 9.04 & 3.05 & \textbf{1.89} \\ 
        mcp124-1$^\divideontimes$ & 124 & 124 & \textbf{$<$0.005} & 0.01 & 2.05 & 0.12 & theta6$^\divideontimes$ & 300 & 4375 & 6.77& 19.49 & 4.22 & \textbf{3.90} \\ 
        mcp124-2$^\divideontimes$ & 124 & 124 & \textbf{$<$0.005} & \textbf{$<$0.005} & 1.90 & 0.13 & theta12$^\star$ & 601 & 17979 & 86.23 & 127  & 11.69 & \textbf{9.85} \\ 
        mcp124-3$^\divideontimes$ & 124 & 124 & \textbf{0.01} & \textbf{0.01} & 1.83 & 0.12 & theta102$^\star$ & 501 & 37467 & 71.44 & 141  & 6.36 & \textbf{4.57} \\ 
        mcp124-4$^\divideontimes$ & 124 & 124 & \textbf{0.01} & \textbf{0.01} & 1.56 & 0.13 & theta123$^\star$ & 601 & 90020 & 206  & 369  & 10.60 & \textbf{8.92} \\ \hline
    \end{tabular}
    \begin{tablenotes}
\footnotesize
\item [1] Problems from SDPLIB are marked with "$^\divideontimes$", and problems from Hans Mittelmann's SDP benchmarks are marked with "$^\star$"
\end{tablenotes}

\end{threeparttable}
\end{table}
}

To further benchmark the performance of \oursolver{}, we employ the scaled shifted geometric means (SSGM) to analyze solving times for problems sourced from both SDPLIB and Hans Mittelmann's benchmarks respectively. The SSGMs calculated for SDPLIB illustrate \oursolver{}'s proficiency in solving small-scale SDP problems, whereas the SSGMs for Hans Mittelmann's benchmarks assess the efficacy of \oursolver{} in addressing larger SDP challenges. For an observation $y \in \mbb {R} ^n_+$, the shifted geometric mean is calculated as follows:
$$
SGM = \exp{\left [\sum_{i=1}^n \frac{1}{n}\log (\max\{1, y_i + s\})\right ]} - s,
$$
where $s$ represents the shift parameter, set to 10 in this case. To compare the solving times of different solvers, we first calculate the shifted geometric mean for each solver's solving time. Then, we derive the scaled shifted geometric mean by normalizing these values against the smallest one. The SSGMs on SDPLIB and Hans Mittelmann's SDP benchmarks are shown in Table \ref{tab:SSGM}.

\ifthenelse{\boolean{Arxiv}}{
\begin{table}[!ht]
    \centering
    \caption{SSGMs of Different Solvers}
    \begin{tabular}{c|cccc}
    \hline
        Problem Set & \oursolver{} & SDPLR & SDPNAL+ & COPT \\ \hline
        SDPLIB & 1.13 & 2.46 & 13.59 & 1 \\ 
         Mittelmann’s Benchmarks & 1 & 3.68 & 10.81 & 1.44 \\ \hline
    \end{tabular}
    \label{tab:SSGM}
\end{table}
}{}

In the evaluation of the solvers' performance on the selected problems from the SDPLIB using SSGMs, COPT emerges as the top performer, closely followed by \oursolver{}. Compared to SDPLR, \oursolver{} demonstrates approximately 1.5 times the speed and outpaces SDPNAL+ by a factor of about 9. For selected problems from Hans Mittelmann’s SDP Benchmarks, \oursolver{} outperforms all competitors, achieving a significant speed advantage over SDPLR by a factor of 3.68 and surpassing SDPNAL+ by 10.81 times. \oursolver{} also marginally exceeds COPT by 1.44 times. These observations suggest that \oursolver{} has a more pronounced advantage in handling larger-scale SDP problems.

\section{Conclusion\label{sec:conclusion}}
In this paper, we have introduced an ADMM approach for solving SDP problems. This method incorporates a matrix variable splitting technique, leading to a low-rank bilinear decomposition problem distinct from traditional \BM{} factorization. We have shown that the ADMM algorithm has a sublinear convergence rate to the KKT solutions and attains a linear rate under the fulfillment of the \L{}ojasiewicz inequality with an exponent of 1/2.
Based on the ADMM algorithm, we developed a novel framework \oursolver{} that employs a \BM{} {\WS} technique alongside the logarithmic rank reduction with a dynamic rank selection strategy ensuring stable performance. Extensive empirical evaluations prove that \oursolver{} can effectively handle various SDP problems and is capable of effectively solving large-scale SDP problems. 

It is important to remark that the development and tuning of \oursolver{} are still at their preliminary level. Future enhancements, such as the adaptive step sizes for the ADMM algorithm, more sophisticated warm-start strategies, intelligent automation for switching phases, and additional heuristic methods, are anticipated to significantly elevate the performance of \oursolver{}. 
Moreover, the current design of \oursolver{} is entirely a first-order framework, which makes it possible for the GPU-based implementation. With the integration of GPU's parallel computing capabilities, it is promising that \oursolver{} may unlock its full potential for addressing super large-scale SDP challenges.

\section*{Acknowledgments}
The authors thank Wenzhi Gao from Stanford University and Mingyu Lei and Jiayu Zhang from Shanghai Jiao Tong University for their insightful discussions and assistance.

\newpage
\bibliographystyle{plainnat}
\bibliography{ref}
\newpage

\appendix
\section{Complete Theoretical Analysis}\label{sec:theory}
\subsection{Proof of Theorem~\ref{prop:KKT}}

\begin{proof}{Proof}
For notation simplicity, we will drop the superscript $\sigma$ in the constants $L_f^\sigma, L^\sigma$ and $L_0^\sigma$ in the proof. 
Firstly, using the $\epsilon$-KKT condition \eqref{def:eKKT-split} of \eqref{prob:BM}, we have
\begin{align}\label{eq:uvdiff}
\begin{split}
        \gamma\|\overline{U}-\overline{V}\|_F \leq & \epsilon+ \left\|\nabla f(\overline{V}\overline{U}^{\top}) \overline{U} +\sum_{i=1}^{m}\overline{\lambda}_iA_{i}\overline{U} \right\|_F \\
        \leq& \epsilon+L_f\norm{\overline{U}}_F + \sqrt{s_A}  \norm{\overline{{\lambda}}}_2\cdot \|\overline{U} \|_F \\
        \leq &\epsilon+ L_f \sigma + \sqrt{s_A}\sigma^2,
        \end{split}
    \end{align}
    where the inequalities holds because $(\overline{U},\overline{V},\overline{\lambda})$ lies in the bounded set $\mcal B_{\sigma}\brbra{0^{n\times r},0^{n\times r},0^{m}}$, $\|\nabla f\|_F$ is bounded above by $L_f$ on $B_{\sigma^2}\brbra{0^{n\times n}}$ and 
    $$
    \left\| \sum_{i=1}^{m}\overline{\lambda}_iA_{i}\overline{U} \right\|_F \leq \sum_{i=1}^{m} |\overline{\lambda}_i| \cdot \|A_{i}\|_2 \cdot \|\overline{U}\|_F \leq \sqrt{s_A}  \norm{\overline{{\lambda}}}_2\cdot \|\overline{U} \|_F 
    $$
    where the first inequality holds because of triangle inequality and the second inequality because of the Cauchy–Schwarz inequality. Next, adding up the first two equations in \eqref{def:eKKT-split} and using the triangle inequality
     yield
    \begin{equation}\label{eq:KKT-dervive}
\Big\|\nabla f(\overline{U}\overline{V}^{\top})\overline{V} + \nabla f(\overline{V}\overline{U}^{\top})\overline{U}  +2 \sum_{i=1}^{m}\overline{\lambda}_iA_{i}\widehat{U}\Big\|_F \leq 2 \epsilon
    \end{equation}
Moreover, it is easily verified that
\begin{align}\label{eq:partialUV}
\begin{split}
   & \nabla f(\overline{U}\overline{V}^{\top})\overline{V} + \nabla f(\overline{V}\overline{U}^{\top})\overline{U}- 2\nabla f(\widehat{U}\widehat{U}^{\top}) \widehat{U} \\
    =&  \left(\nabla f(\overline{U}\overline{V}^{\top})+\nabla f(\overline{V}\overline{U}^{\top}) - 2\nabla f\big(  \frac{\overline{U}\overline{V}^{\top} + \overline{V}\overline{U}^{\top}}{2}
    \big) \right) \widehat{U} \\
    &+ 2\left(\nabla f\big(  {\frac{\overline{U}\overline{V}^{\top} +\overline{V}\overline{U}^{\top}}{2}}
    \big)-\nabla f(\widehat{U}\widehat{U}^{\top})
    \right)\widehat{U} \\
    &+  \nabla f(\overline{U}\overline{V}^{\top})(\overline{V} -\widehat{U} ) + \nabla f(\overline{V}\overline{U}^{\top})(\overline{U} - \widehat{U}) 
\end{split}
\end{align}
Using \cite[Proposition 3.2.12]{ortega2000iterative} and the Lipschitz continuity of $\nabla^2 f$, we know
\[
\|\nabla f(X)-\nabla f(Y)
-\nabla^2 f(Y)(X-Y)\|_F
\leq {L_0\over 2} \|X-Y\|_F^2,
\]
for $X,Y\in \mcal B_{\sigma^2}\brbra{0^{n\times n}}$. Then, by applying this inequality to the pairs $\left(\overline{U}\overline{V}^{\top},\frac{\overline{U}\overline{V}^{\top}+\overline{V}\overline{U}^{\top}}{2}
\right)$ and $\left(\overline{V}\overline{U}^{\top},\frac{\overline{U}\overline{V}^{\top}+\overline{V}\overline{U}^{\top}}{2}\right)$, and then summing up these resulting inequalities, one can easily show that
\begin{align}\label{ineq:part1}
\begin{split}
   &\left\| \nabla f(\overline{U}\overline{V}^{\top}) + \nabla f(\overline{V}\overline{U}^{\top})- 2\nabla f \left(\frac{\overline{U}\overline{V}^{\top} + \overline{V}\overline{U}^{\top}}{2} \right) \right\|_F\\
   \leq\, & \frac{L_0}{4}\|\overline{U}\overline{V}^{\top} - \overline{V}\overline{U}^{\top}\|_F^2\\
   = & \frac{L_0}{4}\|\overline{U}(\overline{V}^{\top}
   -\overline{U}^{\top}) + (\overline{U}
   -\overline{V})\overline{U}^{\top}\|_F^2
   \\
   \leq \,& L_0 \sigma^2\|\overline{U} - \overline{V}\|_F^2,
   \end{split}
\end{align}
where the last inequality holds because $\overline{U}$ and $\overline{V}$ are bounded by $\sigma$. Besides, since $\nabla f$ is Lipschitz continuous, we have
\begin{align}\label{ineq:part2}
\begin{split}
    \left\| \nabla f \left(\frac{\overline{U}\overline{V}^{\top} + \overline{V}\overline{U}^{\top}}{2} \right) -  \nabla f(\widehat{U}\widehat{U}^{\top})\right\|_F &\leq  L \left\| \frac{\overline{U}\overline{V}^{\top} + \overline{V}\overline{U}^{\top}}{2}  - \widehat{U}\widehat{U}^{\top}\right\|_F \\ &=\frac{L}{4}\|(\overline{U} - \overline{V})(\overline{U} - \overline{V})^\top\|_F\\
   &\leq  \frac{L}{4}\|\overline{U} - \overline{V}\|_F^2
    \end{split}
\end{align}
For the third term in~\eqref{eq:partialUV}, we can derive that
\[
\nabla f(\overline{U}\overline{V}^{\top})(\overline{V} -\widehat{U} ) + \nabla f(\overline{V}\overline{U}^{\top})(\overline{U} - \widehat{U}) = \frac{1}{2}\left( \nabla f(\overline{U}\overline{V}^{\top}) - \nabla f(\overline{V}\overline{U}^{\top})\right) (\overline{V} - \overline{U}),
\]
which implies that 
\begin{align}\label{ineq:part3}
\begin{split}
    \norm{\nabla f(\overline{U}\overline{V}^{\top})(\overline{V} -\widehat{U} ) + \nabla f(\overline{V}\overline{U}^{\top})(\overline{U} - \widehat{U})}_F 
    & \leq    \frac{1}{2} \norm{\nabla f(\overline{U}\overline{V}^{\top}) - \nabla f(\overline{V}\overline{U}^{\top})}_F \norm{\overline{V} - \overline{U}}_F  \\
    & \leq \frac{L}{2} \norm{\overline{U}\overline{V}^{\top} - \overline{V}\overline{U}^{\top}}_F \norm{\overline{V} - \overline{U}}_F\\
    & \leq L\sigma \norm{\overline{V} - \overline{U}}_F^2.
    \end{split}
\end{align}
Substituting~\eqref{ineq:part1}, \eqref{ineq:part2} and~\eqref{ineq:part3} into~\eqref{eq:partialUV}, we get 
\begin{align*}
   \| \nabla f(\overline{U}\overline{V}^{\top})\overline{V} + \nabla f(\overline{V}\overline{U}^{\top})\overline{U}- 2\nabla f(\widehat{U}\widehat{U}^{\top}) \widehat{U}\|_F 
  \leq  (L_0 \sigma^3 +\frac{3L\sigma}{2} )\|\overline{U} - \overline{V}\|_F^2,
\end{align*}
which together with~\eqref{eq:uvdiff} and~\eqref{eq:KKT-dervive} yields
\begin{align*}
\Big\|2\nabla f(\widehat{U}\widehat{U}^{\top}) \widehat{U}  +2 \sum_{i=1}^{m}\overline{\lambda}_iA_{i}\widehat{U}\Big\|_F
\leq  2\epsilon + \frac{(2L_0 \sigma^3 +3L \sigma)(\epsilon+ L_f \sigma + \sqrt{s_A}\sigma^2)^2}{2\gamma^2}.
\end{align*}

\noindent Similarly, since $\|{\cal A}(\overline{U}{\overline{V}}^{\top})-b\|_2= \|{\cal A}(\overline{V}{\overline{U}}^{\top}) - b\|_2\leq \epsilon$ and thus it holds
\begin{align*}
\left\| {\cal A}(\widehat{U}{\widehat{U}}^{\top}) - b \right\|_2& \leq \epsilon+{\frac{1}{2}} \left\| 2{\cal A}(\widehat{U}{\widehat{U}}^{\top}) - {\cal A}(\overline{U}{\overline{V}}^{\top})-{\cal A}(\overline{V}{\overline{U}}^{\top}) \right\|_2 \\
&\leq \epsilon+\frac{\|\mathcal{A}\|_2}{2} \cdot \left\|2 \widehat{U}\widehat{U}^{\top} -\overline{U}\overline{V}^{\top}-\overline{V}\overline{U}^{\top} \right\|_F \\
&\leq\epsilon+ \frac{\|
\mathcal{A}\|_2(\epsilon+ L_f \sigma + \sqrt{s_A}\sigma^2)^2}{4\gamma^2}.
\end{align*}
\end{proof}

\subsection{Proof of Proposition~\ref{lem:decent}}

\begin{proof}{Proof}
Note that
\begin{equation*}
\begin{split} & \mcal L_{\rho}(U^{k+1},V^{k+1},{\lambda}^{k+1})-\mcal L_{\rho}(U^{k},V^{k},{\lambda}^{k})\\
= & \underbrace{\mcal L_{\rho}(U^{k+1},V^{k+1},{\lambda}^{k+1})-\mcal L_{\rho}(U^{k+1},V^{k+1},{\lambda}^{k})}_{(\rone)}\\
 & +\underbrace{\mcal L_{\rho}(U^{k+1},V^{k+1},{\lambda}^{k})-\mcal L_{\rho}(U^{k+1},V^{k},{\lambda}^{k})}_{(\rtwo)}\\
 & +\underbrace{\mcal L_{\rho}(U^{k+1},V^{k},{\lambda}^{k})-\mcal L_{\rho}(U^{k},V^{k},{\lambda}^{k})}_{(\rthree)}.
\end{split}
\end{equation*}
Now, we bound term $(\rone)-(\rthree)$ separately,
\begin{equation*}
\begin{split}(\rone)= & \mcal L_{\rho}(U^{k+1},V^{k+1},{\lambda}^{k+1})-\mcal L_{\rho}(U^{k+1},V^{k+1},{\lambda}^{k})\\
= & ({\lambda}^{k+1}-{\lambda}^{k})^{\top}\brbra{{\cal A}(U^{k+1}(V^{k+1})^{\top})-b}\\
= & \rho\norm{{\cal A}(U^{k+1}(V^{k+1})^{\top})-b}_2^{2}\\
= & \frac{1}{\rho}\norm{{\lambda}^{k+1}-{\lambda}^{k}}_2^{2}.
\end{split}
\end{equation*}
Because of the strongly convexity of subproblems, we have
\begin{equation*}
\begin{split}\mcal L_{\rho}(U^{k+1},V^{k},{\lambda}^{k})+\frac{\gamma}{2}\norm{U^{k+1}-U^{k}}_{F}^{2} & \leq\mcal L_{\rho}(U^{k},V^{k},{\lambda}^{k}),\end{split}
\label{eq:threeU}
\end{equation*}
and
\begin{equation*}
\mcal L_{\rho}(U^{k+1},V^{k+1},{\lambda}^{k})+\frac{\gamma}{2}\norm{V^{k+1}-V^{k}}_{F}^{2}\leq\mcal L_{\rho}(U^{k+1},V^{k},{\lambda}^{k}).\label{eq:threeV}
\end{equation*}
We can then bound term $(\rtwo)$ and $(\rthree)$
\begin{equation*}
\begin{split}(\rtwo)= & \mcal L_{\rho}(U^{k+1},V^{k+1},{\lambda}^{k})-\mcal L_{\rho}(U^{k+1},V^{k},{\lambda}^{k})\\
\leq & -\frac{\gamma}{2}\norm{V^{k+1}-V^{k}}_{F}^{2}.
\end{split}
\end{equation*}
Similarly, we have
\begin{equation*}
(\rthree)\leq-\frac{\gamma}{2}\norm{U^{k+1}-U^{k}}_{F}^{2}.
\end{equation*}
By combining the above three inequalities we complete the proof.
\end{proof}

\subsection{Proof of Proposition~\ref{propostion:assu1}}

\begin{proof}{Proof}
For any given $\beta>0$, $S_\beta$ is closed because $f$ is a smooth convex function. Suppose $S_\beta$ is not empty.\\
(i) 
Let $\mu$ be the strongly convex modulus 
of $f(X)+\frac{\rho_0}{2}\|\mcal{A}(X)-b\|_2^2$ and $X^*$ be its minimum point. It follows from the strongly convexity that
\[
f(UV^{\top})+\frac{\rho_0}{2}\|\mcal{A}(UV^{\top})-b\|_2^2
\geq f(X^*) + \frac{\rho_0}{2}\|\mcal{A}(X^*)-b\|_2^2
+\frac{\mu}{2}\|UV^{\top}-X^*\|_F^2,
\]
which implies $\{UV^{\top}\,|\,(U,V)\in S_\beta\}$ is bounded. On the other hand,
$\|U-V\|_F$ is bounded for any $(U,V)\in S_\beta$. Note that
\[
\|U\|_F^2 +\|V\|_F^2=\|U-V\|_F^2 + 2\inprod{U}{V}
\]
We can deduce that $S_\beta$ is bounded. This completes the proof of the first part.\\
(ii) Given the strictly feasible solution $\lambda^*$, we have
\begin{align}\label{eq:cuv}
    \inprod C{UV^\top}
        =\inprod {C-\sum_{i=1}^m y_i^* A_i}{UV^\top}+\sum_{i=1}^my_i^*\inprod {A_i}{UV^\top}
        =\inprod {C-\sum_{i=1}^my_i^* A_i}{UV^\top}+\inprod{{y}^*}{\mcal A(UV^\top)}.
\end{align}
For the first term in \eqref{eq:cuv}, we have
\begin{align*}
   \inprod {C-\sum_{i=1}^m y_i^* A_i}{UV^\top} = & \inprod {C-\sum_{i=1}^m y_i^* A_i}{UU^\top} + \inprod {C-\sum_{i=1}^m y_i^* A_i}{U(V - U)^\top} \\
    \geq & \lambda_{\min}\|U\|_F^2 - \inprod {(C-\sum_{i=1}^m y_i^* A_i)U}{V - U} \\
    \geq & \lambda_{\min}\|U\|_F^2 -\left\|\left(C-\sum_{i=1}^m y_i^* A_i\right)U\right\|_F\|V - U\|_F \\
   \geq & \lambda_{\min}\|U\|_F^2 - \lambda_{\max}\|U\|_F\|V-U\|_F \\
   \geq & \frac{\lambda_{\min}}{2} \|U\|_F^2 - \frac{\lambda_{\max} \kappa_c}{2}\|V-U\|_F^2,
\end{align*}
where the second inequality holds because the fact $\inprod {B_1}{B_2}\geq -\|B_1\|_F\|B_2\|_F$ for any $B_1, B_2$, the third inequality holds because of $\|B_1B_2\|_F\leq \|B_1\|_F \|B_2\|_F$, and the last inequlity holds because $\lambda_{\max}\|U\|_F\|V-U\|_F \leq \frac{\lambda_{\min}}{2} \|U\|_F^2 + \frac{\lambda^2_{\max}}{2\lambda_{\min}}\|V-U\|_F^2$ and $\kappa_c = \frac{\lambda_{\max}}{\lambda_{\min}}$. For the second term in \eqref{eq:cuv}, we have
\begin{align*}
    \inprod{{y}^*}{\mcal A(UV^\top)} =  \inprod{{y}^*}{\mcal A(UV^\top)-b}+\inprod{{y}^*}{b}
        \geq-\frac{1}{2}\norm{{y}^*}_2^2-\frac{1}{2}\norm{{\cal A}(UV^\top)-b}_2^2+\inprod{{y}^*}{b}.
\end{align*}
By choosing $\rho_0 \geq 1 +  \lambda_{\max} \kappa_c$, we have 
\begin{align*}
    & f(UV^{\top})+\frac{\rho_{0}}{2}\norm{U-V}_{F}^{2}+\frac{\rho_{0}}{2}\norm{{\cal A}(UV^{\top})-b}_2^{2} \\
    \geq & \frac{\lambda_{\min}}{2} \|U\|_F^2 +\left(\frac{\rho_0}{2} -  \frac{\lambda_{\max} \kappa_c}{2} \right) \|V-U\|_F^2 + \left(\frac{\rho_0}{2} -  \frac{1}{2} \right) \norm{{\cal A}(UV^\top) -b}_2^2  -\frac{1}{2}\norm{{y}^*}_2^2+ \inprod{{y}^*}{b} \\
    \geq & \frac{\lambda_{\min}}{2} \|U\|_F^2 - \frac{1}{2}\norm{{y}^*}_2^2+ \inprod{{y}^*}{b}.
\end{align*}
  We can deduce from the above inequality that 
\begin{align}\label{ineq:boundu}
    \|U\|_F^2 \leq \frac{1}{\lambda_{\min}}\left( \norm{{y}^*}_2^2 - 2\inprod{{y}^*}{b} + 2\beta \right).
\end{align}
By exchanging $U$ and $V$, it is easy to show that \eqref{ineq:boundu} also holds for $V$. Hence, $S_\beta$ is bounded and thus compact.
\end{proof}

\subsection{Proof of Lemma~\ref{lem:bound of dual variable} }
\begin{proof}{Proof}  We establish part (i) of this lemma by mathematical induction while proving part (ii)-(iv) automatically in the process. 
As $(U^0,V^0,\lambda^0)$ is selected according to \eqref{initial}, part (i) naturally holds for $k=0$. Now, assuming part (i) holds for  $k$-th iterate, we proceed to establish its validity for  $k+1$. By this assumption and Proposition~\ref{lem:decent}, we obtain 
\begin{equation*}
\begin{split} 
\hat{\beta}&>  L_{\rho}(U^{k},V^{k},{\lambda}^{k})
+\frac{1}{\rho}\|\lambda^k\|_2^2
\\
& \geq  L_{\rho}(U^{k+1},V^{k+1},{\lambda}^{k})
+\frac{1}{\rho}\|\lambda^k\|_2^2
\\
& =  f(U^{k+1}(V^{k+1})^{T})+\frac{\gamma}{2}\norm{U^{k+1}-V^{k+1}}_F^{2}+\inprod{{\lambda}^{k}}{\mcal A(U^{k+1}V^{k+1\top})-b}+\frac{\rho}{2}\norm{\mcal A(U^{k+1}V^{k+1\top})-b}_2^{2} +{1\over \rho}\|\lambda^k\|_2^2\\
& \geq  f(U^{k+1}(V^{k+1})^{T})+\frac{\gamma}{2}\norm{U^{k+1}-V^{k+1}}_F^{2}-\frac{\norm{{\lambda}^{k}}_2^{2}}{2(\rho-\rho_{0})}-\frac{(\rho-\rho_{0})}{2}\norm{\mcal A(U^{k+1}V^{k+1\top})-b}_2^{2}+{1\over \rho}\|\lambda^k\|_2^2\\
 &\quad +\frac{\rho}{2}\norm{\mcal A(U^{k+1}V^{k+1\top})-b}_2^{2}\\
 &\geq f(U^{k+1}(V^{k+1})^{T})+\frac{\rho_{0}}{2}\norm{U^{k+1}-V^{k+1}}_F^{2}+\frac{\rho_{0}}{2}\norm{\mcal A(U^{k+1}V^{k+1\top})-b}_2^{2}+ \left({1\over \rho}-{1\over 2(\rho-\rho_{0})}\right)\norm{{\lambda}^{k}}_2^{2}\\
 & \geq f(U^{k+1}(V^{k+1})^{T})+\frac{\rho_{0}}{2}\norm{U^{k+1}-V^{k+1}}_F^{2}+\frac{\rho_{0}}{2}\norm{\mcal A(U^{k+1}V^{k+1\top})-b}_2^{2},
\end{split}
\end{equation*}
where the last inequality follows from $\rho>2\rho_0$. This
implies $(U^{k+1},V^{k+1})\in S_{\hat{\beta}}$
and it thus holds from Assumption \ref{assu:compact} that
$\|U^{k+1}\|_F\leq \delta$ and
$\|V^{k+1}\|_F\leq \delta$. Analogous to the proof above, we can also establish
\[
f(U^{k}(V^{k})^{T})+\frac{\rho_{0}}{2}\norm{U^{k}-V^{k}}_F^{2}+\frac{\rho_{0}}{2}\norm{\mcal A(U^{k}V^{k\top})-b}_2^{2} \leq \hat{\beta}.
\]
It therefore holds $(U^{k+1},V^{k+1})\in S_{\hat{\beta}}$ and $\|U^{k}\|_F\leq \delta,
\|V^{k}\|_F\leq \delta$. 

Next, we bound the dual update $\lambda^k$ and $\lambda^{k+1}-\lambda^k$.
To simplify notation, we will omit the superscript $\delta$ from the constants $L_f^\delta$ and $L^\delta$ throughout the proof. According to the optimality condition \eqref{eq:optV} and the update rule of $\boldsymbol{\lambda}^{k+1}$, we have
\begin{equation}
\nabla f(V^{k+1}(U^{k+1})^{\top})U^{k+1}+\sum_{i=1}^{m}\lambda_{i}^{k+1}A_{i}U^{k+1}-\gamma(U^{k+1}-V^{k+1})=0.\label{eq:16}
\end{equation}
Note that $(U^{k+1},V^{k+1})\in S_{\hat{\beta}}$ and it follows from Assumption \ref{assu:full rank} that
\begin{equation}\label{non-sing}
\sigma_{\min} \left(\mcal C(U^{k+1})\right)\geq \eta.
\end{equation}
Let $B_{k+1}=\nabla f(V^{k+1}(U^{k+1})^{\top})U^{k+1}-\gamma(U^{k+1}-V^{k+1})$,
and $b_{k+1}=\vec(B_{k+1})\in\mbb R^{nr}$. Let $C^{k+1}=\mcal C(U^{k+1})=[\vec{(}A_{1}U^{k+1}),\cdots,\vec{(}A_{m}U^{k+1})]\in\mbb R^{nr\times m}$. Then we can rewrite equation \eqref{eq:16} as
\begin{equation}\label{eq:k+1dual}
C^{k+1}{\lambda}^{k+1}=-b^{k+1}.
\end{equation}
Combining \eqref{non-sing} and \eqref{eq:k+1dual}, one can easily show
\begin{equation}
\begin{split}\eta\norm{{\lambda}^{k+1}}_2\leq&\norm{b^{k+1}}_2= \norm{B^{k+1}}_F\\
\leq & \norm{\nabla f(V^{k+1}(U^{k+1})^{\top})U^{k+1}}_F+\gamma\norm{U^{k+1}-V^{k+1}}_F\\
\leq & L_{f}\delta+2\gamma \delta.
\end{split}
\label{eq:lambda}
\end{equation}
Similarly, we have
\begin{equation*}
C^{k}{\lambda}^{k}=-b^{k}\label{eq:LS_k}
\end{equation*}
and it therefore holds
\begin{equation*}
\begin{split}C^{k+1}{\lambda}^{k+1}-C^{k}{\lambda}^{k} & =C^{k}\Delta{\lambda}^{k+1}+(C^{k+1}-C^{k}){\lambda}^{k+1}=-(b^{k+1}-b^{k}).
\end{split}
\end{equation*}
Then we can bound $\norm{\Delta{\lambda}^{k+1}}_2$ by the following inequality
\begin{equation}
\begin{split}\eta\norm{\Delta{\lambda}^{k+1}}_2\leq \norm{C^{k}\Delta{\lambda}^{k+1}}_2\leq \norm{b^{k+1}-b^{k}}_2+\norm{C^{k+1}-C^{k}}_F\norm{{\lambda}^{k+1}}_2.
\end{split}
\label{eq:delta_lambda}
\end{equation}
In what follows, we bound $\norm{b^{k+1}-b^{k}}_2$ and $\norm{C^{k+1}-C^{k}}_F$
separately
\begin{equation}
\begin{split}\norm{b^{k+1}-b^{k}}_2= & \norm{B^{k+1}-B^{k}}_F\\
\leq & \norm{\nabla f\brbra{V^{k+1}(U^{k+1})^{T}}U^{k+1}-\nabla f\brbra{V^{k}(U^{k})^{T}}U^k}_F+\gamma\norm{\Delta U^{k+1}}_F+\gamma\norm{\Delta V^{k+1}}_F\\
\leq & \brbra{L\delta^{2}+L_{f}+\gamma}\norm{\Delta U^{k+1}}_F+\brbra{L\delta^{2}+\gamma}\norm{\Delta V^{k+1}}_F,
\end{split}
\label{eq:delta_b}
\end{equation}
where the first inequality is because of the definition of $B^k$ and $B^{k+1}$, and the last inequality is due to Assumption \ref{assu:Lipschitz}.
Moreover, according to the definition of $C^{k}$, we have
\begin{equation}
\begin{split}\norm{C^{k+1}-C^{k}}_F^{2} =\sum_{i=1}^{m}\norm{A_{i}\brbra{U^{k+1}-U^{k}}}_F^{2} \leq\sum_{i=1}^{m}\norm{A_{i}}_F^{2}\norm{U^{k+1}-U^{k}}_F^{2} =s_{A}\norm{\Delta U^{k+1}}_F^{2},
\end{split}
\label{eq:delta_C 1}
\end{equation}
which together with \eqref{eq:lambda} yields
\begin{equation}
\begin{split}\norm{C^{k+1}-C^{k}}_F\norm{\lambda^{k}}_2\leq & \frac{1}{\eta}\brbra{L_{f}\delta\sqrt{s_{A}}+2\gamma \delta\sqrt{s_{A}}}\norm{\Delta U^{k+1}}_F\end{split}
\label{eq:C_lambda}
\end{equation}
Combining \eqref{eq:delta_lambda}, \eqref{eq:delta_b} and \eqref{eq:C_lambda}, we have
\begin{equation*}
\begin{split}\eta\norm{{\lambda}^{k+1}-{\lambda}^{k}}_2\leq & \brbra{L\delta^{2}+L_{f}+\gamma+\frac{1}{\eta}\brbra{L_{f}\delta\sqrt{s_{A}}+2\gamma \delta\sqrt{s_{A}}}}\norm{\Delta U^{k+1}}_F\\
 & +\brbra{L\delta^{2}+\gamma}\norm{\Delta V^{k+1}}_F.
\end{split}
\end{equation*}
Therefore, we have 
\begin{equation}
\begin{split}\norm{{\lambda}^{k+1}-{\lambda}^{k}}_2^{2}\leq & \frac{2}{\eta^{2}}\brbra{L\delta^{2}+L_{f}+\gamma+\frac{1}{\eta}L_{f}\delta\sqrt{s_{A}}+\frac{2}{\eta}\gamma \eta\sqrt{s_{A}}}^{2}\norm{\Delta U^{k+1}}_F^{2}\\
 &+\frac{2}{\eta^{2}}\brbra{L\delta^{2}+\gamma}^{2}\norm{\Delta V^{k+1}}_F^{2}\\
\triangleq & c_{1}\norm{\Delta U^{k+1}}_F^{2}+c_{2}\norm{\Delta V^{k+1}}_F^{2}
\end{split}\label{dualbound}
\end{equation}

Finally, we establish the descent of the augmented Lagrangian function and justify the validation of part (i) for iterate $k+1$. By invoking proposition \eqref{lem:decent} and using \eqref{dualbound}, it holds that
\begin{equation*}, 
\begin{split}&\mcal L_{\rho}(U^{k+1},V^{k+1},{\lambda}^{k+1})-\mcal L_{\rho}(U^{k},V^{k},{\lambda}^{k})\\
\leq&  \frac{1}{\rho}\norm{\Delta{\lambda}^{k+1}}_2^{2}-\frac{\gamma}{2}\brbra{\norm{\Delta U^{k+1}}_F^{2}+\norm{\Delta V^{k+1}}^{2}}\\
\leq& -(\frac{\gamma}{2}-\frac{c_{1}}{\rho})\norm{\Delta U^{k+1}}_F^{2}-(\frac{\gamma}{2}-\frac{c_{2}}{\rho})\norm{\Delta V^{k+1}}_F^{2}\\
\leq&  0.
\end{split}
\end{equation*}
Since $\rho> {1\over \hat{\beta}-\beta^0} \left({L_{f}\delta+2\gamma \delta\over \eta}\right)^2$, we know
\begin{equation*}
\begin{split} 
  L_{\rho}(U^{k+1},V^{k+1},{\lambda}^{k+1})
+{1\over \rho}\|\lambda^{k+1}\|_2^2
 &\leq   L_{\rho}(U^{k},V^{k},{\lambda}^{k})
+{1\over \rho}\|\lambda^{k+1}\|_2^2
\\
 &\leq  \beta^0 +{1\over \rho} \left({L_{f}\delta+2\gamma \delta\over \eta}\right)^2
< \hat{\beta}.
\end{split}
\end{equation*}
This completes the proof.

\end{proof}

\subsection{Proof of Theorem~\ref{thm:convergence}}
\begin{proof}{Proof}
Based on the developments outlined in Lemma \ref{lem:bound of dual variable}, we deduce the existence of a constant $\kappa$ such that
\begin{equation}\label{eq:decrease of ALF-revise}
\begin{split}\mcal L_{\rho}(U^{k+1},V^{k+1},{\lambda}^{k+1})-\mcal L_{\rho}(U^{k},V^{k},{\lambda}^{k})\leq -\kappa\norm{\Delta U^{k+1}}_F^{2}-\kappa\norm{\Delta V^{k+1}}_F^{2}.
\end{split}
\end{equation}
On the other hand, we know from
\eqref{eq:optv_revise} and \eqref{eq:u-kkt} that 
\begin{equation*}
\|\nabla \mcal L(U^{k+1},V^{k+1},\lambda^{k+1})\|_F^2 \leq \left({L_f^\delta + L^\delta \delta^2  +  \rho s_A\delta^2+\|\mathcal{A}\|_2\delta }+\gamma\right)^2\|\Delta V^{k+1}\|_F^2 +{1\over \rho^2}\|\Delta \lambda^{k+1}\|^2,
\end{equation*}
where $\mcal L$ is the Lagrangian function of \eqref{prob:penalty}. Note that
\[
\mcal L_\rho(U,V,\lambda)
= \mcal L (U,V,\lambda)+{\rho \over 2}\|\mcal A(UV^T)-b\|_2^2, 
\]
there must exist a constant $\zeta>0$ such that
\begin{equation}\label{gradient}
\|\nabla \mcal L_{\rho}(U^{k+1},V^{k+1},\lambda^{k+1})\|_F   \leq \zeta(\|\Delta U^{k+1}\|_F+\|\Delta V^{k+1}\|_F+\|\Delta \lambda^{k+1}\|_2).
\end{equation}
Therefore, by using \eqref{gradient}, \eqref{eq:decrease of ALF-revise} and Assumption~\ref{assum:KL}, we can invoke \citep[Theorem 2.9]{attouch2013convergence} that
the sequence $\{(U^k,V^k,\lambda^k)\}_{k=0}^{\infty}$ converges to a critical point of $L_\rho(U,V,\lambda)$, which corresponds to a KKT point of problem \eqref{prob:penalty}. We further have the sequence $\{(U^k,V^k,\lambda^k)\}_{k=0}^{\infty}$ has a finite length, i.e.
\[
\sum_{j=1}^\infty \|\Delta U^{j+1}\|_F +\|\Delta V^{j+1}\|_F
+\|\Delta\lambda^{j+1}\|_2< 
\infty.
\]
By elementary calculus, one can show that
\[
\min_{1\leq j\leq K}  \|\Delta U^{j+1}\|_F +\|\Delta V^{j+1}\|_F
+\|\Delta \lambda^{j+1}\|_2 =O({1\over K}).
\]
This, together with \eqref{eq:optv_revise} and \eqref{eq:u-kkt}, implies that
\begin{equation*}
\begin{split}
\min_{1\leq j\leq K} &\left\{\|\nabla f(U^{k+1}(V^{k+1})^{\top})V^{k+1}+\gamma(U^{k+1}-V^{k+1})+\sum_{i=1}^{m}\lambda_{i}^{k+1}A_{i}V^{k+1}\|_F\right.\\
&\,+\|\nabla f(V^{k+1}(U^{k+1})^{\top})V^{k+1}+\gamma(V^{k+1}-U^{k+1})+\sum_{i=1}^{m}\lambda_{i}^{k+1}A_{i}V^{k+1}\|_F \\
&\,+\|\mcal A(U^{k+1}(V^{k+1})^T)-b\|_2\bigg\} =O\left(\frac{1}{K}\right).
\end{split}
\end{equation*}
In other words, Algorithm~\ref{alg:ADMM-for-SDP} requires $O(1/\epsilon)$ iterations to produce an $\epsilon$-KKT point of \eqref{prob:penalty}.
\end{proof}

\subsection{Proof of Theorem~\ref{thm:linear convergence}}
\begin{proof}{Proof}
Because of the convergence of the sequence $\{(U^k,V^k,\lambda^k)\}_{k=0}^{\infty}$, the iteration point will enter the local area around the limit point and satisfies the \L{}ojasiewicz inequality after a constant number of steps. For the simplicity of notation, we assume $(U^k,V^k,\lambda^k)\in \mcal B_{\delta_c}(\overline{U},\overline{V},\overline{\lambda})$ for any $k\geq 0$ in the following analysis. 

First, consider the case where $\alpha =0$. Assume that the sequence does not converge in finite steps and hence $\mcal L_{\rho}(U^k, V^k,\lambda^k)-\mcal L_{\rho}(\overline{U},\overline{V},\overline{\lambda})>0$ holds for any $k\geq 0$. According to the \L{}ojasiewicz inequality, we then have
\begin{align*}
    c\norm{\nabla\mcal L_{\rho}(U^k,V^k,\lambda^k)}_F\geq1.
\end{align*}
This contradicts to the result that the sequence converges to a KKT point. Hence, the sequence must converge in finite steps when $\alpha=0$.

Next, we consider the remaining cases where $\alpha\in(0,1)$. Define $\phi(t)=\frac{c}{1-\alpha}t^{1-\alpha}$. We can then rewrite the \L{}ojasiewicz inequality \eqref{eq:KL inequality} as
\begin{equation}\label{eq:KL phi}
    \phi'(\mcal L_{\rho}(U,V,{\lambda})-\mcal L_{\rho}(\overline{U},\overline{V},\overline{\lambda}))\cdot\norm{\nabla \mcal L_{\rho}(U,V,{\lambda})}_F\geq1.
\end{equation}
Let $l^k=\mcal L_{\rho}(U^{k},V^{k},{\lambda}^{k})-\mcal L_{\rho}(\overline{U},\overline{V},\overline{\lambda})$. Then according to the \L{}ojasiewicz inequality, we have
\begin{align*}
    \phi'(l^k)\cdot\norm{\nabla \mcal L_{\rho}(U^k,V^k,{\lambda}^k)}_F\geq 1.
\end{align*}
Because $l^k-l^{k+1}\geq0$, by multiplying it to both side we get the following inequality
\begin{align*}
    \phi'(l^k)(l^k-l^{k+1})\cdot\norm{\nabla \mcal L_{\rho}(U^k,V^k,{\lambda}^k)}_F\geq l^k-l^{k+1}.
\end{align*}
Note that $\phi$ is concave, we further have
\begin{align*}
    \brbra{\phi(l^k)-\phi(l^{k+1})}\cdot\norm{\nabla \mcal L_{\rho}(U^k,V^k,{\lambda}^k)}_F\geq l^k-l^{k+1}.
\end{align*}
According to Lemma~\ref{lem:bound of dual variable} and~\eqref{gradient},  we can bound $l^k-l^{k+1}$ and $\norm{\nabla \mcal L_{\rho}(U^k,V^k,{\lambda}^k)}_F$ respectively
\begin{align*}
    l^k-l^{k+1}\geq \kappa\brbra{\norm{\Delta U^{k+1}}_F^2+\norm{\Delta V^{k+1}}_F^2},\\
    \norm{\nabla \mcal L_{\rho}(U^k,V^k,{\lambda}^k)}_F\leq \zeta\brbra{\norm{\Delta U^{k}}_F+\norm{\Delta V^{k}}_F+\|\Delta \lambda^{k}\|_2}.
\end{align*}
Note that $\norm{\Delta{\lambda}^{k+1}}_2^{2}$ can be upper bounded by $\norm{\Delta U^{k+1}}_F^{2}$ and $\norm{\Delta V^{k+1}}_F^{2}$, hence there exists some constant $\tilde{\kappa}>0$ such that 
\begin{align*}
     l^k-l^{k+1}\geq \tilde{\kappa}\brbra{\norm{\Delta U^{k+1}}_F^2+\norm{\Delta V^{k+1}}_F^2+\|\Delta \lambda^{k+1}\|_2^2}.
\end{align*}
Combining the above inequalities, the following inequality follows naturally
\begin{align*}
    \frac{\zeta}{\tilde\kappa}\brbra{\phi(l^k)-\phi(l^{k+1})}\brbra{\norm{\Delta U^{k}}_F+\norm{\Delta V^{k}}_F+\|\Delta \lambda^{k}\|_2}\geq \brbra{\norm{\Delta U^{k+1}}_F^2+\norm{\Delta V^{k+1}}_F^2+\|\Delta \lambda^{k+1}\|_2^2}.
\end{align*}
Recall the inequality that $2\sqrt{ab}\leq a+b$, we can further get
\begin{align*}
    \frac{2\zeta}{\tilde\kappa}\brbra{\phi(l^k)-\phi(l^{k+1})}+\brbra{\norm{\Delta U^{k}}_F+\norm{\Delta V^{k}}_F+\|\Delta \lambda^{k}\|_2}\geq 2\brbra{\norm{\Delta U^{k+1}}_F+\norm{\Delta V^{k+1}}_F+\|\Delta \lambda^{k+1}\|_2}.
\end{align*}
By adding the above inequality from $k$ to $K$, we have
\begin{align*}
    \frac{2\zeta}{\tilde\kappa}\brbra{\phi(l^{k})-\phi(l^{K+1})}+\brbra{\norm{\Delta U^{k}}_F+\norm{\Delta V^{k}}_F+\|\Delta \lambda^{k}\|_2}\geq& \sum_{t=k}^{K}\brbra{\norm{\Delta U^{t+1}}_F+\norm{\Delta V^{t+1}}_F+\|\Delta \lambda^{t+1}\|_2}\\
    &+\brbra{\norm{\Delta U^{K+1}}_F+\norm{\Delta V^{K+1}}_F+\|\Delta \lambda^{K+1}\|_2}.
\end{align*}
Because $\mcal L_{\rho}(U^{k},V^{k},{\lambda}^{k})$ is non-increasing, we have $l^{K+1}>0$ and hence $\phi(l^{k})-\phi(l^{K+1})\leq \phi(l^{k})$. Taking $K$ to $\infty$, we can derive the following inequality
\begin{align}
    \frac{2\zeta}{\tilde\kappa}\phi(l^{k})+\brbra{\norm{\Delta U^{k}}_F+\norm{\Delta V^{k}}_F+\|\Delta \lambda^{k}\|_2}\geq \sum_{t=k}^{\infty}\brbra{\norm{\Delta U^{t+1}}_F+\norm{\Delta V^{t+1}}_F+\|\Delta \lambda^{t+1}\|_2}.
\end{align}
Let $S^{k}=\sum_{t=k}^{\infty} \brbra{\norm{\Delta U^{t+1}}_F+\norm{\Delta V^{t+1}}_F+\|\Delta \lambda^{t+1}\|_2}$. Then we can rewrite the above inequality as follows
\begin{align}\label{eq:sum of length}
    \frac{2\zeta}{\tilde\kappa}\phi(l^{k})+S^{k-1}-S^{k}\geq S^{k}.
\end{align}
According to the definition of $\phi$ and the \L{}ojasiewicz inequality, we have
\begin{align*}
    c(l^k)^{-\alpha}\cdot\norm{\nabla \mcal L_{\rho}(U^k,V^k,{\lambda}^k)}_F\geq 1,
\end{align*}
and hence there exists a constant $\tilde c>0$ such that
\begin{align*}
    \phi(l^k)=\frac{c}{1-\alpha}(l^k)^{1-\alpha}\leq \tilde c\brbra{\norm{\Delta U^{k}}_F+\norm{\Delta V^{k}}_F+\|\Delta \lambda^{k}\|_2}^{\frac{1-\alpha}{\alpha}}=\tilde c(S^{k-1}-S^k)^{\frac{1-\alpha}{\alpha}}.
\end{align*}
Combining with \eqref{eq:sum of length}, we have
\begin{align*}
    \frac{2\zeta\tilde{c} }{\tilde\kappa}(S^{k-1}-S^k)^{\frac{1-\alpha}{\alpha}}+S^{k-1}-S^{k}\geq S^{k}.
\end{align*}
The convergence property of such sequence has been studied in Theorem 2 of \cite{attouch2009convergence}, which gives the analysis of convergence rate.

\end{proof}

\section{Additional Numerical Results}\label{app:B}
In this section, we provide the solving errors for the instances tested in Section~\ref{sec:numerical}.
Since each solver terminates with distinct criteria, we present three commonly defined errors for the solutions of each solver including:

Primal Infeasibility:
$$\frac{\norm{\mcal A(X)-b}_2}{1+\norm{b}_{1}},$$

Dual Infeasibility: 
$$\frac{\abs{\min\{0,\sigma_{\min}(C-\mcal A^{*}(\lambda))\}}}{1+\norm c_{1}},$$

Primal Dual Gap:
$$\frac{\inner{C}{X}-\lambda^{\top}b}{1+\abs{\inner{C}{X}}+\abs{\lambda^{\top}b}}.$$ 

The three error measures are chosen for ranking in Hans Mittelmann’s benchmarks, and it is also widely used in many existing solvers, like COPT \citep{ge2022cardinal} and HDSDP~\citep{gao2022hdsdp}. The reports of errors in all three experiments are shown in the following tables.

\ifthenelse{\boolean{Arxiv}}{\begin{table}[!ht]
    \centering
    \scriptsize
    \caption{Solving Errors on Max-cut Problems}
    \setlength{\tabcolsep}{0.85mm}
    \begin{tabular}{c|ccc|ccc|ccc|ccc}
    \hline
        ~ & ~ & \oursolver{} & ~ & ~ & SDPLR & ~ & ~ & SDPNAL+ & ~ & ~ & COPT & ~ \\ \hline
         Problem & P Infeasi & D Infeasi & PD Gap & P Infeasi & D Infeasi & PD Gap & P Infeasi & D Infeasi & PD Gap & P Infeasi & D Infeasi & PD Gap \\ 
        G1 & 1.19E-08 & 2.51E-07 & 4.10E-09 & 1.84E-08 & 4.51E-09 & 1.59E-06 & 8.11E-08 & 2.27E-08 & 2.64E-07 & 6.54E-13 & 0.00E+00 & 2.37E-09 \\ 
        G2 & 1.18E-08 & 2.36E-07 & 4.90E-09 & 2.44E-08 & 6.09E-09 & 1.77E-06 & 1.53E-07 & 2.76E-08 & 1.35E-06 & 1.60E-13 & 0.00E+00 & 3.21E-09 \\ 
        G3 & 1.16E-08 & 5.62E-07 & 3.21E-09 & 1.28E-08 & 2.94E-08 & 3.81E-07 & 2.59E-07 & 3.03E-08 & 4.90E-06 & 1.00E-13 & 0.00E+00 & 3.93E-09 \\ 
        G4 & 1.35E-08 & 4.58E-07 & 3.68E-09 & 6.88E-09 & 1.54E-08 & 4.12E-07 & 1.95E-07 & 1.91E-08 & 1.37E-06 & 1.85E-13 & 0.00E+00 & 2.25E-09 \\ 
        G5 & 1.33E-08 & 3.21E-07 & 6.70E-09 & 2.29E-08 & 7.51E-09 & 2.58E-07 & 2.77E-07 & 8.80E-08 & 6.39E-06 & 5.84E-14 & 0.00E+00 & 3.38E-09 \\ 
        G6 & 5.97E-09 & 4.57E-07 & 8.81E-09 & 1.62E-08 & 8.56E-09 & 9.76E-08 & 6.95E-08 & 7.13E-09 & 9.55E-07 & 5.20E-14 & 0.00E+00 & 2.39E-08 \\ 
        G7 & 8.77E-09 & 1.62E-07 & 2.17E-08 & 2.46E-08 & 5.77E-10 & 6.35E-06 & 6.43E-08 & 6.53E-09 & 1.06E-06 & 4.00E-12 & 0.00E+00 & 1.22E-08 \\ 
        G8 & 6.89E-09 & 6.61E-07 & 1.29E-08 & 2.03E-08 & 1.02E-08 & 3.64E-06 & 5.13E-08 & 2.87E-09 & 5.39E-07 & 7.10E-13 & 0.00E+00 & 1.11E-08 \\ 
        G9 & 5.21E-09 & 1.00E-07 & 9.17E-09 & 1.77E-08 & 4.50E-09 & 1.48E-06 & 2.05E-09 & 3.51E-08 & 1.38E-05 & 2.60E-13 & 0.00E+00 & 1.15E-08 \\ 
        G10 & 3.94E-09 & 8.37E-08 & 8.66E-09 & 1.44E-08 & 5.20E-09 & 1.92E-06 & 6.09E-08 & 2.00E-08 & 1.04E-05 & 1.80E-13 & 0.00E+00 & 9.39E-09 \\ 
        G11 & 1.58E-08 & 2.95E-07 & 2.38E-08 & 2.25E-08 & 1.04E-08 & 9.26E-07 & 2.53E-08 & 5.69E-08 & 1.19E-05 & 1.53E-13 & 0.00E+00 & 6.44E-08 \\ 
        G12 & 1.88E-08 & 4.21E-07 & 3.48E-08 & 2.40E-08 & 1.56E-08 & 9.59E-07 & 2.23E-07 & 5.84E-08 & 1.61E-05 & 1.34E-13 & 0.00E+00 & 2.78E-08 \\ 
        G13 & 1.87E-08 & 1.64E-07 & 3.11E-08 & 2.04E-08 & 1.43E-08 & 3.31E-06 & 1.04E-15 & 7.66E-08 & 1.26E-05 & 1.90E-14 & 0.00E+00 & 4.09E-08 \\ 
        G14 & 2.46E-08 & 2.00E-06 & 1.73E-07 & 1.61E-08 & 8.89E-09 & 3.72E-07 & 4.27E-08 & 8.15E-08 & 2.83E-06 & 2.29E-14 & 0.00E+00 & 2.45E-08 \\ 
        G15 & 1.17E-08 & 1.24E-06 & 9.35E-08 & 2.11E-08 & 1.42E-08 & 1.36E-06 & 1.22E-07 & 2.59E-08 & 8.42E-07 & 1.67E-12 & 0.00E+00 & 6.14E-09 \\ 
        G16 & 6.75E-09 & 1.00E-06 & 1.08E-08 & 2.18E-08 & 1.23E-08 & 1.24E-06 & 1.61E-07 & 1.02E-07 & 7.55E-06 & 1.51E-13 & 0.00E+00 & 1.01E-08 \\ 
        G17 & 1.72E-08 & 1.17E-06 & 9.45E-08 & 2.32E-08 & 2.11E-08 & 3.54E-06 & 5.46E-08 & 2.09E-08 & 8.83E-07 & 2.87E-13 & 0.00E+00 & 7.50E-09 \\ 
        G18 & 1.80E-08 & 2.95E-07 & 1.41E-07 & 1.62E-08 & 1.16E-08 & 1.96E-06 & 3.56E-07 & 2.51E-09 & 3.78E-06 & 6.69E-13 & 0.00E+00 & 3.36E-08 \\ 
        G19 & 3.04E-09 & 2.55E-07 & 3.37E-09 & 1.81E-08 & 2.27E-08 & 4.63E-06 & 4.27E-08 & 7.12E-10 & 1.98E-07 & 4.40E-13 & 0.00E+00 & 6.31E-08 \\ 
        G20 & 2.49E-08 & 3.67E-07 & 1.45E-07 & 1.75E-08 & 2.58E-07 & 1.48E-06 & 1.03E-07 & 1.53E-08 & 2.73E-07 & 3.54E-13 & 0.00E+00 & 4.74E-08 \\ 
        G21 & 2.32E-08 & 1.26E-07 & 1.73E-07 & 1.97E-08 & 6.01E-09 & 5.85E-07 & 2.28E-08 & 5.29E-08 & 8.87E-06 & 7.33E-13 & 0.00E+00 & 4.10E-08 \\ 
        G22 & 2.33E-09 & 3.68E-07 & 2.58E-09 & 7.46E-09 & 9.26E-09 & 1.12E-06 & 1.93E-07 & 4.47E-08 & 9.03E-06 & 3.24E-14 & 0.00E+00 & 4.80E-09 \\ 
        G23 & 2.09E-09 & 4.76E-07 & 1.07E-09 & 5.03E-09 & 4.28E-09 & 6.25E-07 & 1.27E-08 & 3.63E-08 & 9.80E-06 & 3.28E-14 & 0.00E+00 & 7.53E-09 \\ 
        G24 & 2.35E-09 & 2.29E-07 & 2.96E-09 & 9.46E-09 & 7.39E-09 & 1.25E-06 & 1.90E-08 & 6.64E-08 & 1.28E-05 & 1.73E-13 & 0.00E+00 & 1.79E-09 \\ 
        G25 & 2.47E-09 & 2.95E-07 & 2.95E-09 & 6.63E-09 & 8.38E-09 & 9.15E-07 & 1.88E-07 & 5.96E-08 & 1.04E-05 & 3.29E-14 & 0.00E+00 & 5.88E-09 \\ 
        G26 & 2.39E-09 & 2.87E-07 & 3.27E-09 & 9.27E-09 & 1.78E-09 & 1.28E-06 & 2.04E-07 & 4.47E-08 & 1.00E-05 & 6.83E-15 & 0.00E+00 & 8.11E-09 \\ 
        G27 & 8.27E-09 & 9.04E-07 & 1.29E-08 & 6.69E-09 & 1.49E-09 & 9.32E-07 & 7.83E-08 & 3.89E-09 & 1.83E-06 & 8.81E-14 & 0.00E+00 & 7.47E-09 \\ 
        G28 & 8.76E-09 & 7.36E-07 & 5.41E-08 & 5.30E-09 & 6.75E-09 & 3.86E-07 & 7.45E-08 & 1.38E-08 & 5.95E-06 & 2.19E-14 & 0.00E+00 & 4.63E-09 \\ 
        G29 & 9.94E-10 & 6.00E-07 & 1.80E-09 & 4.24E-09 & 0.00E+00 & 3.25E-06 & 4.18E-08 & 2.55E-08 & 9.31E-06 & 2.38E-13 & 0.00E+00 & 4.80E-09 \\ 
        G30 & 1.11E-09 & 6.00E-07 & 1.96E-09 & 8.81E-09 & 1.15E-08 & 9.01E-07 & 6.79E-08 & 3.16E-09 & 1.57E-06 & 2.48E-12 & 0.00E+00 & 9.40E-09 \\ 
        G31 & 7.24E-09 & 1.03E-06 & 2.80E-08 & 6.52E-09 & 8.20E-09 & 9.62E-07 & 3.02E-16 & 2.51E-08 & 1.05E-05 & 8.41E-14 & 0.00E+00 & 7.05E-09 \\ 
        G32 & 1.07E-09 & 2.86E-08 & 2.95E-10 & 8.51E-09 & 6.59E-10 & 1.94E-06 & 7.41E-08 & 3.47E-08 & 2.03E-05 & 2.80E-14 & 0.00E+00 & 4.15E-08 \\ 
        G33 & 1.04E-09 & 2.84E-08 & 3.43E-09 & 9.34E-09 & 3.66E-09 & 1.98E-07 & 1.34E-15 & 3.69E-08 & 1.90E-05 & 6.40E-14 & 0.00E+00 & 4.54E-08 \\ 
        G34 & 9.31E-10 & 1.63E-08 & 1.95E-09 & 5.07E-09 & 1.10E-09 & 8.99E-07 & 6.08E-09 & 3.81E-08 & 1.58E-05 & 2.79E-14 & 0.00E+00 & 1.10E-08 \\ 
        G35 & 9.53E-09 & 6.24E-07 & 9.16E-08 & 8.67E-09 & 6.81E-09 & 2.24E-06 & 4.90E-16 & 7.34E-08 & 2.72E-05 & 1.49E-13 & 0.00E+00 & 6.46E-09 \\ 
        G36 & 4.29E-09 & 5.41E-07 & 5.29E-08 & 8.23E-09 & 4.77E-09 & 5.69E-07 & 9.81E-09 & 2.73E-08 & 1.11E-05 & 1.78E-13 & 0.00E+00 & 7.59E-09 \\ 
        G37 & 8.49E-09 & 7.48E-07 & 8.77E-08 & 8.63E-09 & 0.00E+00 & 6.44E-06 & 8.40E-08 & 3.55E-08 & 4.64E-06 & 8.52E-13 & 0.00E+00 & 5.63E-09 \\ 
        G38 & 8.52E-09 & 4.37E-07 & 1.25E-07 & 9.13E-09 & 3.75E-09 & 3.85E-07 & 4.50E-08 & 4.98E-08 & 1.08E-05 & 6.95E-13 & 0.00E+00 & 4.93E-09 \\ 
        G39 & 1.74E-09 & 5.79E-08 & 8.31E-09 & 9.05E-09 & 2.03E-09 & 1.63E-06 & 2.45E-08 & 4.39E-08 & 1.45E-05 & 6.48E-13 & 0.00E+00 & 6.79E-09 \\ 
        G40 & 3.25E-09 & 1.13E-07 & 3.04E-08 & 6.87E-09 & 6.19E-09 & 5.75E-07 & 3.59E-16 & 4.55E-08 & 1.07E-05 & 9.76E-13 & 0.00E+00 & 6.76E-09 \\ 
        G41 & 3.03E-09 & 1.93E-07 & 3.23E-08 & 5.55E-09 & 9.54E-09 & 1.62E-06 & 1.35E-07 & 1.37E-08 & 4.93E-06 & 1.86E-14 & 0.00E+00 & 2.88E-08 \\ 
        G42 & 2.71E-09 & 5.15E-08 & 2.03E-08 & 9.77E-09 & 1.93E-08 & 1.70E-06 & 2.80E-16 & 4.38E-08 & 1.62E-05 & 7.22E-15 & 0.00E+00 & 7.04E-09 \\ 
        G43 & 5.24E-09 & 2.87E-07 & 7.00E-09 & 1.60E-08 & 5.26E-09 & 1.01E-06 & 3.76E-16 & 5.98E-06 & 1.45E-03 & 8.51E-14 & 0.00E+00 & 8.73E-09 \\ 
        G44 & 6.01E-09 & 2.94E-07 & 7.80E-10 & 1.38E-08 & 3.04E-09 & 5.46E-07 & 3.03E-16 & 6.18E-06 & 1.38E-03 & 8.90E-14 & 0.00E+00 & 8.15E-09 \\ 
        G45 & 3.93E-09 & 3.14E-07 & 2.50E-09 & 1.67E-08 & 4.54E-09 & 1.40E-06 & 1.36E-08 & 7.65E-09 & 3.66E-07 & 1.71E-13 & 0.00E+00 & 6.09E-09 \\ 
        G46 & 4.53E-09 & 2.43E-07 & 4.14E-09 & 7.87E-09 & 3.52E-09 & 4.46E-07 & 2.77E-07 & 2.23E-09 & 1.65E-06 & 3.92E-12 & 0.00E+00 & 7.03E-09 \\ 
        G47 & 5.61E-09 & 2.71E-07 & 4.12E-09 & 1.10E-08 & 1.35E-08 & 4.29E-07 & 1.18E-07 & 5.26E-09 & 1.35E-07 & 8.55E-13 & 0.00E+00 & 3.92E-09 \\ 
        G48 & 2.72E-09 & 5.76E-09 & 1.48E-08 & 5.90E-09 & 1.39E-10 & 1.88E-07 & 1.26E-14 & 3.36E-07 & 3.95E-03 & 1.65E-17 & 0.00E+00 & 5.28E-09 \\ 
        G49 & 3.33E-09 & 5.77E-09 & 1.30E-08 & 6.11E-09 & 3.26E-10 & 5.35E-07 & 6.03E-14 & 3.28E-07 & 5.26E-04 & 1.71E-17 & 0.00E+00 & 5.28E-09 \\ 
        G50 & 2.38E-09 & 5.57E-08 & 1.54E-08 & 5.24E-09 & 0.00E+00 & 2.26E-07 & 8.06E-15 & 2.28E-07 & 1.85E-04 & 6.43E-17 & 0.00E+00 & 1.26E-08 \\ 
        G51 & 1.83E-08 & 1.01E-06 & 7.32E-08 & 1.24E-08 & 1.86E-08 & 6.79E-07 & 2.14E-07 & 3.54E-08 & 7.63E-07 & 1.30E-13 & 0.00E+00 & 9.00E-09 \\ 
        G52 & 1.81E-08 & 1.17E-06 & 1.37E-07 & 1.80E-08 & 2.19E-08 & 1.69E-06 & 6.27E-08 & 2.05E-08 & 6.22E-07 & 2.90E-14 & 0.00E+00 & 9.60E-09 \\ 
        G53 & 1.62E-08 & 5.05E-07 & 7.00E-09 & 1.58E-08 & 1.34E-08 & 1.02E-06 & 1.28E-15 & 7.01E-06 & 1.61E-03 & 2.77E-14 & 0.00E+00 & 4.73E-09 \\ 
        G54 & 1.35E-08 & 1.05E-06 & 6.46E-08 & 1.58E-08 & 0.00E+00 & 4.91E-06 & 1.60E-07 & 5.27E-08 & 1.66E-06 & 5.65E-13 & 0.00E+00 & 5.02E-09 \\ 
        G55 & 1.66E-09 & 1.35E-07 & 1.08E-08 & 3.12E-09 & 1.04E-09 & 4.29E-07 & 1.81E-08 & 2.66E-08 & 1.44E-06 & 9.48E-15 & 0.00E+00 & 7.23E-09 \\ 
        G56 & 3.65E-09 & 2.49E-07 & 2.94E-08 & 1.76E-09 & 1.71E-09 & 8.54E-07 & 4.88E-09 & 1.31E-08 & 7.14E-06 & 6.56E-15 & 0.00E+00 & 5.04E-09 \\ 
        G57 & 9.60E-10 & 9.40E-09 & 4.69E-09 & 3.12E-09 & 8.73E-10 & 3.24E-07 & 9.90E-10 & 2.02E-08 & 2.79E-05 & 2.19E-13 & 0.00E+00 & 1.12E-08 \\ 
        G58 & 1.76E-09 & 3.31E-07 & 2.14E-08 & 1.97E-09 & 9.71E-11 & 3.70E-06 & 5.60E-08 & 4.08E-08 & 4.05E-05 & 7.96E-14 & 0.00E+00 & 6.47E-09 \\ 
        G59 & 1.92E-09 & 1.21E-07 & 2.87E-08 & 3.74E-09 & 7.40E-10 & 1.04E-06 & 2.66E-08 & 1.40E-08 & 1.28E-05 & 2.79E-13 & 0.00E+00 & 8.42E-09 \\ 
        G60 & 9.24E-10 & 1.63E-07 & 4.42E-09 & 2.00E-09 & 5.05E-10 & 6.38E-07 & t & t & t & 1.25E-15 & 0.00E+00 & 4.15E-09 \\ 
        G61 & 1.22E-09 & 2.26E-07 & 4.63E-09 & 2.79E-09 & 6.91E-10 & 2.16E-06 & t & t & t & 1.16E-15 & 0.00E+00 & 8.51E-09 \\ 
        G62 & 7.77E-10 & 3.72E-09 & 4.36E-08 & 1.14E-09 & 2.15E-10 & 1.42E-07 & t & t & t & 4.89E-14 & 0.00E+00 & 6.76E-09 \\ 
        G63 & 8.78E-10 & 2.02E-07 & 5.25E-09 & 2.52E-09 & 1.26E-09 & 3.67E-07 & t & t & t & 1.45E-14 & 0.00E+00 & 4.03E-09 \\ 
        G64 & 2.19E-09 & 1.92E-07 & 4.28E-08 & 1.75E-09 & 3.87E-09 & 6.00E-07 & t & t & t & 4.54E-13 & 0.00E+00 & 2.99E-09 \\ 
        G65 & 5.79E-10 & 3.74E-09 & 9.13E-09 & 1.40E-09 & 3.05E-11 & 5.35E-07 & t & t & t & 6.44E-14 & 0.00E+00 & 5.59E-09 \\ 
        G66 & 4.45E-10 & 2.76E-09 & 1.00E-08 & 3.37E-10 & 2.06E-10 & 5.27E-08 & t & t & t & 5.78E-14 & 0.00E+00 & 9.57E-09 \\ 
        G67 & 6.68E-10 & 2.98E-09 & 2.04E-08 & 1.19E-09 & 0.00E+00 & 9.29E-07 & t & t & t & 1.20E-13 & 0.00E+00 & 8.31E-09 \\ \hline
    \end{tabular}
    \label{tab:max-cut-err}
\end{table}

\begin{table}[!ht]
    \centering
    \scriptsize
    \caption{Solving Errors on Matrix Completion Problems}
    \setlength{\tabcolsep}{0.80mm}
    \begin{tabular}{c|ccc|ccc|ccc|ccc}
    \hline
        ~ & ~ & \oursolver{} & ~ & ~ & SDPLR & ~ & ~ & SDPNAL+ & ~ & ~ & COPT & ~ \\ \hline
        Problem & P Infeasi & D Infeasi & PD Gap & P Infeasi & D Infeasi & PD Gap & P Infeasi & D Infeasi & PD Gap & P Infeasi & D Infeasi & PD Gap \\ 
        MC\_1000 & 1.45E-10 & 4.04E-08 & 3.84E-06 & 3.00E-10 & 3.54E-07 & 2.26E-05 & 1.71E-18 & 3.69E-09 & 1.49E-06 & 1.30E-06 & 0.00E+00 & 1.96E-08 \\ 
        MC\_2000 & 1.45E-10 & 2.18E-08 & 1.90E-06 & 1.20E-10 & 2.28E-08 & 2.39E-06 & 1.87E-09 & 2.25E-09 & 2.00E-06 & 2.02E-06 & 0.00E+00 & 6.99E-09 \\ 
        MC\_3000 & 5.49E-11 & 4.22E-10 & 9.30E-07 & 4.17E-11 & 3.68E-08 & 1.60E-05 & 1.30E-09 & 3.88E-10 & 1.45E-06 & 4.76E-06 & 0.00E+00 & 4.85E-09 \\ 
        MC\_4000 & 4.58E-11 & 2.54E-09 & 5.78E-06 & 5.40E-11 & 6.14E-08 & 4.07E-06 & 1.62E-09 & 1.23E-09 & 2.98E-06 & 6.77E-06 & 0.00E+00 & 3.11E-09 \\ 
        MC\_5000 & 3.18E-11 & 1.46E-09 & 1.24E-06 & 2.07E-11 & 4.39E-08 & 1.88E-06 & 6.60E-10 & 2.81E-10 & 2.11E-07 & 6.73E-06 & 0.00E+00 & 1.34E-08 \\ 
        MC\_6000 & 3.74E-11 & 3.89E-09 & 2.24E-06 & 4.19E-11 & 1.51E-08 & 1.05E-05 & 2.87E-10 & 0.00E+00 & 1.81E-06 & 5.06E-06 & 0.00E+00 & 9.80E-09 \\ 
        MC\_8000 & 2.57E-11 & 3.75E-09 & 7.26E-07 & t & t & t & 8.84E-11 & 1.29E-08 & 6.39E-06 & t & t & t \\ 
        MC\_10000 & 2.13E-11 & 3.05E-09 & 1.20E-06 & t & t & t & 2.79E-09 & 3.49E-11 & 1.37E-05 & t & t & t \\ 
        MC\_12000 & 1.76E-11 & 4.15E-09 & 2.25E-06 & - & - & - & - & - & - & - & - & - \\ 
        MC\_14000 & 1.26E-11 & 4.52E-10 & 1.48E-06 & - & - & - & - & - & - & - & - & - \\ 
        MC\_16000 & 1.20E-11 & 5.25E-09 & 5.23E-06 & - & - & - & - & - & - & - & - & - \\ 
        MC\_18000 & 1.17E-11 & 2.06E-09 & 4.33E-06 & - & - & - & - & - & - & - & - & - \\ 
        MC\_20000 & 1.09E-11 & 1.53E-09 & 9.54E-07 & - & - & - & - & - & - & - & - & - \\ 
        MC\_40000 & 5.17E-12 & 8.41E-11 & 1.37E-06 & - & - & - & - & - & - & - & - & - \\ \hline
    \end{tabular}
    \label{tab:mc-error}
\end{table}

\begin{table}[!ht]
    \centering
    \scriptsize
    \caption{Solving Errors on Problems from SDPLIB and Mittelmann's Benchmark}
    \setlength{\tabcolsep}{0.70mm}
    \begin{tabular}{c|ccc|ccc|ccc|ccc}
    \hline
        ~ & ~ & \oursolver{} & ~ & ~ & SDPLR & ~ & ~ & SDPNAL+ & ~ & ~ & COPT & ~ \\ \hline
        Problem & P Infeasi & D Infeasi & PD Gap & P Infeasi & D Infeasi & PD Gap & P Infeasi & D Infeasi & PD Gap & P Infeasi & D Infeasi & PD Gap \\ 
        ALH & 1.29E-06 & 1.46E-04 & 1.23E-03 & 1.36E-06 & 2.85E-05 & 4.87E-04 & 7.24E-07 & 1.85E-07 & 6.70E-06 & 7.30E-10 & 0.00E+00 & 2.51E-09 \\ 
        BH2 & 6.07E-07 & 1.82E-02 & 5.30E-03 & 1.04E-06 & 5.07E-05 & 1.24E-04 & 2.72E-07 & 9.39E-08 & 1.78E-06 & 6.21E-10 & 0.00E+00 & 2.54E-09 \\ 
        cancer\_100 & 2.86E-09 & 3.28E-04 & 1.77E-05 & t & t & t & 2.61E-08 & 9.24E-10 & 4.43E-05 & 1.97E-05 & 0.00E+00 & 4.54E-09 \\ 
        CH2 & 8.95E-07 & 2.09E-03 & 6.09E-04 & 1.13E-06 & 1.99E-04 & 2.49E-04 & 7.92E-09 & 2.60E-07 & 3.48E-04 & 6.36E-10 & 0.00E+00 & 3.56E-09 \\ 
        checker\_1.5 & 4.58E-09 & 5.29E-08 & 1.54E-07 & 4.11E-09 & 1.83E-08 & 1.36E-04 & 9.12E-08 & 2.27E-07 & 1.25E-03 & 2.70E-13 & 0.00E+00 & 4.53E-08 \\ 
        G40\_mb & 9.72E-09 & 1.29E-07 & 9.23E-09 & 9.89E-09 & 8.57E-09 & 1.07E-06 & 9.06E-10 & 0.00E+00 & 3.33E-07 & 6.01E-07 & 0.00E+00 & 1.07E-05 \\ 
        G48\_mb & 4.03E-09 & 8.63E-10 & 1.37E-07 & 4.05E-09 & 0.00E+00 & 2.96E-04 & 1.14E-11 & 1.00E-07 & 7.15E-02 & 1.92E-04 & 0.00E+00 & 4.85E-03 \\ 
        gpp100 & 1.96E-07 & 1.41E-06 & 1.75E-06 & 1.98E-07 & 1.56E-08 & 2.62E-06 & 1.65E-12 & 5.20E-06 & 7.97E-04 & 1.67E-08 & 0.00E+00 & 6.59E-09 \\ 
        gpp124-1 & 1.58E-07 & 3.02E-06 & 8.02E-06 & 1.60E-07 & 5.62E-06 & 2.81E-03 & 1.10E-06 & 0.00E+00 & 3.47E-06 & 3.46E-07 & 0.00E+00 & 4.77E-08 \\ 
        gpp124-2 & 1.59E-07 & 8.72E-07 & 1.54E-06 & 1.54E-07 & 2.81E-08 & 1.34E-06 & 9.19E-07 & 0.00E+00 & 4.22E-06 & 4.20E-08 & 0.00E+00 & 2.09E-09 \\ 
        gpp124-3 & 1.58E-07 & 2.52E-06 & 6.50E-07 & 1.57E-07 & 1.41E-08 & 6.04E-07 & 3.32E-12 & 3.12E-06 & 4.59E-04 & 4.96E-08 & 0.00E+00 & 6.00E-09 \\ 
        gpp124-4 & 1.60E-07 & 2.44E-07 & 9.99E-07 & 1.58E-07 & 1.11E-06 & 4.81E-05 & 6.48E-09 & 0.00E+00 & 3.73E-06 & 6.50E-08 & 0.00E+00 & 3.19E-09 \\ 
        gpp250-1 & 7.87E-08 & 4.75E-07 & 2.64E-06 & 7.41E-08 & 1.11E-06 & 3.17E-05 & 6.32E-07 & 0.00E+00 & 9.06E-06 & 7.89E-08 & 0.00E+00 & 2.93E-08 \\ 
        gpp250-2 & 7.75E-08 & 2.90E-07 & 1.20E-06 & 7.42E-08 & 1.78E-09 & 7.79E-07 & 1.79E-11 & 1.05E-05 & 4.93E-03 & 4.53E-08 & 0.00E+00 & 7.70E-09 \\ 
        gpp250-3 & 7.84E-08 & 9.26E-08 & 3.28E-07 & 7.97E-08 & 1.73E-07 & 4.45E-05 & 6.35E-07 & 0.00E+00 & 6.11E-06 & 7.52E-08 & 0.00E+00 & 5.34E-09 \\ 
        gpp250-4 & 7.89E-08 & 5.64E-08 & 3.48E-07 & 7.92E-08 & 6.03E-08 & 1.37E-05 & 3.30E-12 & 8.06E-06 & 1.28E-03 & 1.21E-07 & 0.00E+00 & 3.61E-09 \\ 
        gpp500-1 & 3.98E-08 & 2.00E-07 & 7.84E-07 & 3.88E-08 & 2.37E-07 & 2.11E-04 & 8.89E-07 & 1.15E-11 & 5.76E-06 & 1.22E-06 & 0.00E+00 & 1.02E-05 \\ 
        gpp500-2 & 3.97E-08 & 1.08E-07 & 5.30E-07 & 3.94E-08 & 1.02E-09 & 3.64E-07 & 1.51E-06 & 7.64E-11 & 5.08E-07 & 1.37E-07 & 0.00E+00 & 1.30E-06 \\ 
        gpp500-3 & 3.93E-08 & 3.45E-07 & 2.24E-08 & 3.98E-08 & 2.38E-08 & 1.18E-05 & 1.29E-10 & 5.03E-06 & 4.10E-03 & 1.02E-07 & 0.00E+00 & 1.01E-08 \\ 
        gpp500-4 & 3.98E-08 & 1.17E-07 & 3.41E-08 & 3.99E-08 & 0.00E+00 & 2.34E-05 & 4.75E-13 & 3.95E-10 & 6.99E-07 & 3.23E-07 & 0.00E+00 & 2.47E-06 \\ 
        H3O & 1.13E-06 & 1.29E-04 & 2.47E-04 & 1.18E-06 & 2.13E-05 & 7.14E-05 & 6.84E-08 & 9.82E-08 & 5.39E-06 & 7.36E-10 & 0.00E+00 & 3.49E-09 \\ 
        hand & 7.62E-08 & 2.27E-09 & 2.71E-08 & 1.54E-08 & 4.34E-09 & 1.08E-07 & 8.87E-09 & 3.97E-12 & 2.14E-07 & 2.57E-06 & 0.00E+00 & 4.09E-05 \\ 
        ice\_2.0 & 1.41E-09 & 1.63E-07 & 2.60E-06 & 5.17E-10 & 1.63E-07 & 1.01E-05 & t & t & t & 1.02E-11 & 0.00E+00 & 1.98E-08 \\ 
        mcp100 & 1.03E-07 & 5.71E-06 & 2.60E-07 & 5.96E-08 & 7.03E-08 & 4.04E-06 & 4.90E-16 & 8.26E-06 & 3.55E-05 & 1.13E-11 & 0.00E+00 & 1.72E-09 \\ 
        mcp124-1 & 5.04E-08 & 9.37E-06 & 2.52E-07 & 9.96E-08 & 3.21E-08 & 1.47E-05 & 2.80E-07 & 6.27E-08 & 8.10E-07 & 1.42E-11 & 0.00E+00 & 5.65E-10 \\ 
        mcp124-2 & 1.47E-07 & 1.03E-05 & 3.15E-07 & 1.56E-07 & 2.43E-07 & 1.12E-06 & 4.28E-16 & 1.45E-05 & 3.17E-04 & 5.39E-12 & 0.00E+00 & 6.00E-10 \\ 
        mcp124-3 & 1.55E-07 & 4.64E-06 & 7.67E-08 & 1.44E-07 & 3.99E-08 & 4.80E-06 & 3.66E-16 & 9.71E-06 & 3.94E-04 & 1.74E-12 & 0.00E+00 & 4.61E-09 \\ 
        mcp124-4 & 1.46E-07 & 1.94E-06 & 6.70E-08 & 3.49E-08 & 8.71E-08 & 2.79E-07 & 1.88E-15 & 1.02E-05 & 5.01E-04 & 1.44E-11 & 0.00E+00 & 7.56E-10 \\ 
        mcp250-1 & 4.91E-08 & 2.62E-06 & 2.01E-07 & 3.89E-08 & 4.57E-08 & 5.62E-06 & 2.30E-07 & 1.68E-08 & 1.10E-07 & 3.11E-12 & 0.00E+00 & 4.52E-09 \\ 
        mcp250-2 & 6.07E-08 & 5.33E-06 & 1.53E-07 & 7.23E-08 & 2.70E-09 & 1.29E-05 & 2.01E-15 & 9.92E-06 & 1.09E-04 & 3.97E-12 & 0.00E+00 & 1.26E-09 \\ 
        mcp250-3 & 7.96E-08 & 1.59E-06 & 6.79E-08 & 7.63E-08 & 3.10E-08 & 3.09E-06 & 3.92E-07 & 3.76E-10 & 9.89E-08 & 1.72E-10 & 0.00E+00 & 4.77E-09 \\ 
        mcp250-4 & 3.11E-08 & 1.49E-06 & 3.13E-08 & 7.20E-08 & 8.55E-08 & 4.10E-06 & 3.65E-07 & 1.99E-10 & 1.89E-07 & 1.03E-11 & 0.00E+00 & 9.63E-10 \\ 
        mcp500-1 & 2.71E-08 & 1.63E-06 & 1.16E-07 & 3.25E-08 & 8.11E-08 & 4.40E-06 & 1.25E-07 & 9.61E-08 & 4.60E-06 & 1.19E-15 & 0.00E+00 & 2.23E-07 \\ 
        mcp500-2 & 2.37E-09 & 5.15E-07 & 8.19E-09 & 2.88E-08 & 1.25E-07 & 7.21E-06 & 3.52E-14 & 1.00E-07 & 2.80E-06 & 4.12E-14 & 0.00E+00 & 5.99E-08 \\ 
        mcp500-3 & 6.45E-09 & 1.01E-06 & 1.18E-08 & 2.41E-08 & 6.80E-08 & 8.77E-06 & 2.75E-07 & 2.41E-08 & 1.61E-06 & 9.74E-14 & 0.00E+00 & 2.28E-08 \\ 
        mcp500-4 & 1.97E-08 & 6.44E-07 & 2.43E-08 & 2.54E-08 & 1.18E-08 & 1.29E-06 & 3.10E-16 & 8.22E-06 & 1.04E-03 & 2.00E-14 & 0.00E+00 & 1.10E-08 \\ 
        NH2 & 1.26E-06 & 3.10E-02 & 1.68E-02 & 1.17E-06 & 5.52E-05 & 8.48E-05 & 3.85E-15 & 2.11E-06 & 9.10E-06 & 1.84E-09 & 0.00E+00 & 2.51E-09 \\ 
        NH3 & 8.59E-07 & 1.60E-04 & 9.01E-05 & 1.07E-06 & 3.96E-05 & 2.42E-04 & 1.55E-08 & 1.05E-07 & 2.65E-05 & 5.06E-10 & 0.00E+00 & 2.98E-09 \\ 
        p\_auss2\_3.0 & 1.40E-09 & 1.35E-09 & 9.42E-07 & 1.80E-09 & 1.98E-09 & 3.14E-05 & t & t & t & 1.31E-13 & 0.00E+00 & 4.74E-06 \\ 
        qap5 & 2.53E-06 & 1.07E-03 & 4.60E-05 & 2.44E-06 & 0.00E+00 & 3.08E-03 & 5.52E-07 & 0.00E+00 & 7.87E-09 & 7.91E-08 & 0.00E+00 & 1.89E-09 \\ 
        qap6 & 2.51E-06 & 2.34E-04 & 1.38E-03 & 2.52E-06 & 1.32E-07 & 9.92E-04 & 9.87E-06 & 1.34E-06 & 7.16E-04 & 4.10E-08 & 0.00E+00 & 1.08E-07 \\ 
        qap7 & 2.49E-06 & 1.47E-04 & 2.77E-04 & 2.50E-06 & 6.07E-08 & 9.73E-04 & 8.58E-06 & 4.14E-07 & 5.16E-04 & 3.56E-08 & 0.00E+00 & 6.62E-08 \\ 
        qap8 & 2.48E-06 & 3.24E-04 & 3.59E-04 & 2.49E-06 & 4.89E-08 & 1.32E-03 & 2.96E-06 & 1.15E-07 & 3.29E-04 & 4.88E-08 & 0.00E+00 & 4.96E-08 \\ 
        qap9 & 2.39E-06 & 1.08E-04 & 7.38E-03 & 2.48E-06 & 1.82E-08 & 9.18E-04 & 3.48E-06 & 4.03E-08 & 1.76E-04 & 1.49E-07 & 0.00E+00 & 2.91E-08 \\ 
        qap10 & 2.18E-06 & 1.85E-04 & 1.57E-02 & 2.48E-06 & 2.24E-07 & 1.51E-03 & 6.65E-06 & 2.42E-08 & 3.41E-04 & 1.05E-07 & 0.00E+00 & 3.48E-08 \\ 
        sensor\_500 & 1.92E-07 & 1.02E-03 & 2.06E-04 & 1.93E-07 & 3.22E-04 & 1.77E-06 & 2.94E-07 & 3.50E-09 & 7.60E-06 & 5.85E-08 & 0.00E+00 & 4.41E-10 \\ 
        shmup4 & 2.16E-08 & 1.00E-06 & 6.89E-03 & 2.46E-08 & 5.16E-08 & 2.19E-05 & 2.36E-08 & 3.03E-06 & 3.19E-01 & 3.59E-09 & 0.00E+00 & 4.48E-06 \\ 
        theta1 & 9.96E-06 & 1.51E-04 & 6.33E-06 & 9.99E-06 & 8.87E-07 & 1.59E-05 & 7.71E-06 & 4.39E-08 & 5.09E-06 & 2.83E-08 & 0.00E+00 & 1.60E-09 \\ 
        theta2 & 1.00E-05 & 1.64E-04 & 6.94E-06 & 9.97E-06 & 2.48E-08 & 2.23E-05 & 2.02E-13 & 2.25E-07 & 9.09E-06 & 7.41E-09 & 0.00E+00 & 4.11E-10 \\ 
        theta3 & 9.97E-06 & 5.24E-04 & 6.30E-06 & 9.95E-06 & 5.08E-09 & 2.98E-05 & 2.30E-06 & 2.82E-08 & 1.46E-06 & 2.88E-08 & 0.00E+00 & 9.03E-10 \\ 
        theta4 & 1.00E-05 & 2.57E-04 & 9.27E-07 & 9.94E-06 & 1.94E-09 & 2.35E-05 & 6.13E-13 & 1.54E-08 & 1.12E-04 & 1.01E-08 & 0.00E+00 & 3.56E-10 \\ 
        theta5 & 9.99E-06 & 2.65E-04 & 5.10E-06 & 9.99E-06 & 3.32E-09 & 2.97E-05 & 3.20E-13 & 5.96E-08 & 1.06E-03 & 7.99E-09 & 0.00E+00 & 2.74E-10 \\ 
        theta6 & 9.99E-06 & 4.24E-04 & 9.51E-06 & 9.96E-06 & 1.11E-09 & 2.73E-05 & 6.65E-13 & 4.83E-08 & 7.14E-04 & 7.71E-09 & 0.00E+00 & 2.78E-10 \\ 
        theta12 & 3.09E-05 & 1.80E-03 & 4.12E-05 & 1.00E-05 & 1.39E-09 & 1.83E-05 & 5.02E-06 & 1.15E-09 & 5.31E-06 & 1.74E-08 & 0.00E+00 & 2.12E-09 \\ 
        theta102 & 1.08E-04 & 3.02E-03 & 6.78E-04 & 9.98E-06 & 4.08E-09 & 2.53E-05 & 1.12E-06 & 7.03E-10 & 9.80E-08 & 2.80E-09 & 0.00E+00 & 2.70E-09 \\ 
        theta123 & 2.24E-05 & 8.26E-04 & 2.04E-05 & 9.95E-06 & 2.22E-09 & 3.50E-05 & 3.86E-07 & 6.21E-09 & 5.29E-08 & 1.90E-08 & 0.00E+00 & 2.78E-09 \\ \hline
    \end{tabular}
    \label{tab:sdplib-hans-err}
\end{table}
}{}

\end{document}